\documentclass[11pt]{imsart}

\RequirePackage{amsthm,amsmath,amsfonts,amssymb}
\RequirePackage[numbers]{natbib}

\usepackage{hyperref,diagrams}

\theoremstyle{plain}
\newtheorem{theorem}{Theorem}[section]

\newtheorem{thm}[theorem]{Theorem}
\newtheorem{prop}[theorem]{Proposition}
\newtheorem{lemma}[theorem]{Lemma}
\newtheorem{cor}[theorem]{Corollary}

\theoremstyle{remark}
\newtheorem{defn}[theorem]{Definition}
\newtheorem{remark}[theorem]{Remark}
\newtheorem{notation}[theorem]{Notation}

\renewcommand{\theequation}{\arabic{section}.\arabic{equation}}

\renewcommand{\Bar}{\overline}
\newcommand{\iso}{\cong}
\newcommand{\caniso}{\mbox{\hspace{1ex}}\iso_{\mbox{\hspace{-5ex}\phantom{$\int$}
       \scriptsize{canon.}}}}
\newcommand{\bfr}{{\bf R}}
\newcommand{\disjoint}{{\mbox{\small $\coprod$}}} %disjoint union
\newcommand{\dotprod}{{\mbox{
\,\hspace{-.06in}\raisebox{-.015in}{\bf {\Large $\cdot$}}}}}
\newcommand{\id}{{\rm id.}}
\newcommand{\intersect}{\cap}
\newcommand{\Intersect}{\bigcap}
\newcommand{\plus}{\oplus}
\newcommand{\Plus}{\bigoplus}
\newcommand{\Sym}{{\rm Sym}}
\newcommand{\supp}{{\rm supp}}
\renewcommand{\Tilde}{\widetilde}
\newcommand{\tensor}{\otimes}
\newcommand{\union}{\mbox{\small $\bigcup$}}
\newcommand{\be}{\begin{equation}}
\newcommand{\ee}{\end{equation}}
\newcommand{\bearray}{\begin{eqnarray}}
\newcommand{\eearray}{\end{eqnarray}}
\newcommand{\bestar}{\begin{eqnarray*}}
\newcommand{\eestar}{\end{eqnarray*}}
\newcommand{\ben}{\begin{displaymath}}
\newcommand{\een}{\end{displaymath}}

\newcommand{\qedns}{{\mbox{\hspace{1ex}}
\mbox{\hfill\vrule height10pt width10pt} \\ \vspace{0.5 ex}}}

\newcommand{\pf}{\noi{\bf Proof}: }

\renewcommand{\a}{\alpha}
\renewcommand{\b}{\beta}
\renewcommand{\d}{\delta}
\newcommand{\e}{\epsilon}
\newcommand{\g}{\gamma}
\newcommand{\w}{\omega}
\newcommand{\downto}{\downarrow}
\newcommand{\grad}{\mbox{\rm grad\phantom{-}}}
\newcommand{\hess}{{\rm Hess}}
\newcommand{\na}{\nabla}
\newcommand{\noi}{\noindent}
\renewcommand{\ss}{\vspace{.1in}}
\newcommand{\bs}{\vspace{.2in}}

\newcommand{\bsn}{\vspace{.2in}\noindent}

\newcommand{\ssn}{\vspace{.1in}\noindent}
\newcommand{\End}{{\rm End}}
\newcommand{\lnbrk}{\linebreak}

\newcommand{\rinj}{{r_{\rm inj}}}
\newcommand{\rcx}{{r_{\rm cx}}}
\newcommand{\Hor}{{\rm Hor}}
\newcommand{\Ver}{{\rm Vert}}
\newcommand{\hor}{{\rm hor}}
\newcommand{\ver}{{\rm vert}}
\newcommand{\plast}{{\pi_{\rm last}}}
\newcommand{\tgam}{\tilde{\g}}
\newcommand{\tg}{\tilde{g}}
\newcommand{\F}{{\mathcal{F}}}
\newcommand{\FM}{{\rm FM}}
\newcommand{\EFM}{{\rm EFM}}
\newcommand{\con}{{\rm con}}
\newcommand{\SR}{{\textsc {sr}}}
\newcommand{\PSR}{{\textsc {psr}}}
\newcommand{\dsr}{d_\SR}
\newcommand{\dsrtop}{\dsr^{\rm top}}
\newcommand{\dtop}{d^{\rm top}}
\newcommand{\dpsr}{d_\PSR}

\newcommand{\quo}{{\rm quo}}
\newcommand{\bgp}{\b_m}
\newcommand{\G}{{\mathcal{G}}}

\newcommand{\bgam}{\bar{\g}}
\newcommand{\bF}{\bar{\F}}
\newcommand{\trho}{\tilde{\rho}}

\newcommand{\cS}{S^M}
\newcommand{\cX}{X^M}
\newcommand{\sptop}{S_m^{\rm top}}
\newcommand{\splow}{S_m^{\rm low}}
\newcommand{\sympp}{{\rm Sym^+}(m)}
\newcommand{\devt}{d_{\rm evt}}
\newcommand{\bFtop}{\bF_{\rm top}}
\newcommand{\D}{{\mathcal D}}
\newcommand{\Dtan}{\D^{\rm tan}}
\newcommand{\tD}{\tilde{\D}}
\newcommand{\tDtan}{\tD^{\rm tan}}
\newcommand{\CL}{{\mathcal{C}}}
\newcommand{\CLtan}{\CL^{\rm tan}}

\newcommand{\tmn}{(TM)^{\times N}}
\newcommand{\tf}{\tilde{f}}
\newcommand{\tF}{\tilde{F}}
\newcommand{\ts}{\tilde{s}}
\newcommand{\tz}{\tilde{Z}}
\newcommand{\pb}{{\mathcal{PB}}}
\newcommand{\ord}{{\rm ord}}
\newcommand{\sing}{{\rm sing}}
\newcommand{\clord}{\CL_\ord}
\newcommand{\clsing}{\CL_\sing}

\newcommand{\clsingtan}{\CL_\sing^{\tan}}
\newcommand{\V}{{\mathcal V}}
\newcommand{\B}{{\mathcal B}}
\newcommand{\tB}{\tilde{\B}}
\newcommand{\td}{{\tilde{d}}}
\newcommand{\tE}{{\tilde{E}}}
\newcommand{\tM}{\Tilde{\M}}
\newcommand{\M}{{\mathcal M}}
\newcommand{\N}{{\mathcal N}}
\newcommand{\tim}{{\Tilde{M}}}
\newcommand{\timtop}{\tim^{\rm top}}
\newcommand{\timlow}{\tim^{\rm low}}

\newcommand{\tp}{{\tilde{p}}}
\newcommand{\tq}{{\tilde{q}}}
\newcommand{\tQ}{{\tilde{Q}}}
\newcommand{\tS}{{\tilde{S}}}
\newcommand{\tv}{{\tilde{v}}}
\newcommand{\tX}{{\tilde{X}}}
\newcommand{\Diag}{{\rm Diag}}
\newcommand{\call}{{\mathcal L}}
\newcommand{\cals}{{\mathcal S}}

\newcommand{\fsrs}{f^{\SR}_{\cals}}
\newcommand{\fpsrs}{\tilde{f}^{\PSR}_{\cals}}
\newcommand{\minus}{\,\backslash\, }

\DeclareMathOperator*{\argmin}{arg\,min}

\begin{document}

\begin{frontmatter}

\title{A genericity property of Fr\'{e}chet sample means on Riemannian manifolds
%\\ {\tiny version date: 9/11/2023 } 
}
\runtitle{Genericity of sample means on manifolds}

\author{\fnms{David} \snm{Groisser}\corref{}\ead[label=e1]{groisser@ufl.edu}}

\address{Department of Mathematics\\
University of Florida\\
Gainesville, FL 32611, USA\\
\printead{e1}}
\affiliation{University of Florida}

\and

\author{\fnms{Sungkyu} \snm{Jung}\ead[label=e2]{sungkyu@snu.ac.kr}}

\address{Department of Statistics\\ Seoul National University
\\ Seoul, Korea
\\
\printead{e2}}
\affiliation{Seoul National University}

\and

\author{\fnms{Armin} \snm{ Schwartzman}\ead[label=e4]{armins@ucsd.edu}}
\address{
{Division of Biostatistics}
\\ {and Hal{\i}c{\i}o\u{g}lu Data Science Institute}
\\ {University of California, San Diego}
\\ {La Jolla, CA 92903, USA}
\\ \printead{e4}}
\affiliation{University of California, San Diego}

\runauthor{Groisser, Jung, and Schwartzman}

\begin{abstract}  Let $(M,g)$ be a Riemannian manifold.  If $\mu$ is a probability measure on $M$
given by a continuous density function, one would expect the Fr\'{e}chet means of data-samples 
$Q=(q_1,q_2,\dots, q_N)\in M^N$, with respect to $\mu$, to behave ``generically''; 
e.g. the probability that the Fr\'{e}chet mean set $\FM(Q)$ has
any elements that lie in a given, positive-codimension submanifold, should be zero
for any $N\geq 1$. Even this simplest instance of genericity does not
seem to have been proven in the literature, except in special cases. 
The main result of this paper is a general, and stronger, genericity property:
given i.i.d. absolutely continuous $M$-valued random variables $X_1,\dots, X_N$,
and a subset 
$A\subset M$ of volume-measure zero, 
${\rm Pr}\left\{\FM(\{X_1,\dots,X_N\})\subset M\backslash A\right\}=1.$
We also establish a companion theorem for {\em equivariant Fr\'{e}chet means},
defined when $(M,g)$ arises as  the quotient of a Riemannian manifold $(\tim,\tg)$ by a free,
isometric action of a finite group.  The equivariant Fr\'{e}chet means 
lie in $\tim$, but, as we show, project down to the ordinary Fr\'{e}chet sample means, 
and enjoy a similar genericity property.  Both these theorems are proven as consequences of a purely geometric 
(and quite general) result that constitutes the core mathematics in this paper:  
If $A\subset M$ has volume zero in $M$ , then the set $\{Q\in M^N : \FM(Q) \intersect A\neq\emptyset\}$ 
has volume zero in $M^N$.  We conclude the paper with an application 
to {\em partial scaling-rotation means}, a type of mean for symmetric positive-definite matrices.

\end{abstract}

\begin{keyword}[class=MSC2020]
\kwd[Primary ]{62R30}
%62R30 Statistics on manifolds
%60D05 Geometric probability and stochastic geometry
\kwd{60D05}
\kwd[; secondary ]{60B05,57R99}
%60B05 Probability measures on topological spaces
%57R99 None of the above, but in this section (Differential topology, 57Rxx)
\end{keyword}


\end{frontmatter}

\section{Introduction}
\label{sect:intro}

Let $(M,g)$ be a smooth ($C^\infty$) Riemannian manifold, with induced distance-function $d$, and for
$p\in M$ let $r_p=d(\cdot,p)$ denote the function ``distance to $p$''.   The
{\em Fr\'{e}chet mean set} $\FM(\mu)$ of a Borel probability measure $\mu$ on $M$ is defined to
be the set of minimizers of the function on $M$ defined by $p\mapsto 
\int_M r_p^2\, d\mu$, provided that the integral is finite for some (hence all) $p$, a condition met automatically by all measures $\mu_Q$ considered below.
If $(M,g)$ is complete,  which we will henceforth assume, then $\FM(\mu)$ is nonempty 
for any such $\mu$ \cite[Theorem 2.1]{bhatpat1}.
 Elements of $\FM(\mu)$ are called {\em Fr\'{e}chet means}, or
{\em Riemannian centers of mass}, of $\mu$.  
(In this paper, the term {\em barycenter} will be used
for something more general, defined later.)

For $q\in M$ let $\d_q$ be the Dirac measure with support $\{q\}$, and for $N\geq 1$ let
$M^N$ denote the $N$-fold Cartesian product of $M$ with itself.
For an $N$-tuple  (or {\em configuration}) of points $Q=(q_1,\dots, q_N)\in M^N$, by
the {\em Fr\'{e}chet mean set of $Q$} we will mean $\FM(Q):=\FM(\mu_Q)$, where $\mu_Q=\frac{1}{N}
\sum_{i=1}^N\d_{q_i}$ (the notation $\mu_Q$ will be used with this meaning throughout
the paper); similarly, a {\em Fr\'{e}chet mean
of $Q$} will mean an element of $\FM(Q)$.  For random $N$-tuples $Q$ representing i.i.d.
samples of an underlying probability measure $\mu$, the Fr\'{e}chet means $\mu_Q$ are
called {\em Fr\'{e}chet sample means} of $\mu$.

If $\mu$ is given by a continuous probability density function, one would expect the Fr\'{e}chet
sample means of $\mu$ to behave ``generically''; e.g. the probability that they lie in a
positive-codimension submanifold should be zero
for any $N\geq 1$. Surprisingly, even this much genericity does not
seem to have been proven in the literature, except in special cases.  The main result of this paper is a general, and stronger, genericity property:

\begin{thm}\label{thm1.1} Let $X_1,\dots, X_N$ be i.i.d. 
absolutely continuous $M$-valued random variables.
If
$A\subset M$ has volume zero in $M$, then
\be\label{eq:thm1.1}
{\rm Pr}\left\{\FM(\{X_1,\dots,X_N\})\subset M\backslash A\right\}=1.
\ee
In particular, \eqref{eq:thm1.1} holds if $A$ is a countable union of submanifolds
of positive codimension.
\end{thm}

\noi 
(Terminology in this theorem is discussed the end of this section.)

Fr\'{e}chet means are not always unique (i.e. the Fr\'{e}chet mean set may have more than one element).
In \cite[Theorem 2.1]{afsari}, Afsari established what is, to date, the least-stringent 
 ``concentration'' condition known to guarantee uniqueness.
 Combining this with Theorem
\ref{thm1.1}, we obtain a  corresponding sufficient condition for $N$ i.i.d. random variables
under which, with probability 1, the Fr\'{e}chet 
sample mean 
is unique and lies in the complement
of a given volume-zero subset of $M$. Theorem \ref{thm1.1} and this uniqueness
result are restated together as Corollary \ref{cor:measzero.3}.

 Theorem \ref{thm1.1} has an ``equivariant'' companion, 
in which $(M,g)$ arises as  the quotient of a Riemannian manifold $(\tim,\tg)$ by a free,
isometric action of a finite group $G$. (Thus $(\tim,\tilde{g})$ is a finite-degree  principal Riemannian
cover of $(M,g)$; see Remark \ref{princ_rc}.)  In this setting we consider a hybrid ``equivariant distance
function'' $\devt: M\times \tim\to\bfr$, defined by $\devt(q,\tilde{p})
=\tilde{d}(\quo^{-1}(q),\tilde{p})$, the distance in $(\tim,\tg)$ between the point $\tilde{p}\in\tim$ and the $G$-orbit $\quo^{-1}(q)\subset\tim$, where $\quo:\tim\to M$ is the quotient map.  (Equivalently,
$\devt(q,\tilde{p})=d(q,\quo(\tilde{p}))$.) For a given $Q\in M^N$, 
 we can lift the Fr\'{e}chet objective function $f_Q=\frac{1}{N}\sum_{i=1}^N d(\cdot,q_i)^2 :M\to\bfr$
to a function $\tf_Q:\tim\to\bfr$ by  replacing $d(\cdot, p)$ by $\devt(\cdot, \tilde{p})$.
We call minimizers of such a function {\em equivariant Fr\'{e}chet means} of $Q$
(a special case of ``$\rho$-means'' as extended by Huckemann \cite{huck3DProcrustes2011,Huckemann2011consistency};
if  $\tilde{p}$ is a minimizer, then so is $h\dotprod\tilde{p}$ for every $h\in G$.
For $Q\in M^N$, let ${\rm EFM}(Q)\subset\tim$ denote the set of equivariant Fr\'{e}chet means of $Q$.
We prove the following  result about equivariant Fr\'{e}chet means of
i.i.d. random variables  (essentially an equivariant version of Corollary \ref{cor:measzero.3}, mentioned above):

\begin{thm}
\label{thm1.2} 
Assume that the finite group $G$ acts freely and isometrically on the Riemannian
manifold $(\tim, \tilde{g})$, with quotient $(M,g)$  and $\quo:\tim\to M$
the quotient map.
Let $X, X_1, X_2,\dots, X_N$ be i.i.d. $M$-valued random variables,
with induced probability measure $\mu_X$\,.
 Then:

\begin{itemize}
\item[(a)] $\EFM(\{X_1,\dots,X_N\})=\quo^{-1}(\FM(\{X_1,\dots,X_N\}).$

\ss \item[(b)]  Let $\tilde{A}\subset\tim$ be a $G$-invariant set of 
volume zero in $\tim$, let $A=\quo(\tilde{A})$
and write $\tim_*=\tim\setminus\tilde{A}$ and $M_*=\quo(\tim_*)=M\setminus A.$
If 
$X$ is absolutely continuous then, almost surely, 
 $\FM(\{X_1,\dots,X_N\})\subset M_*$ and $\EFM(\{X_1,\dots,X_N\})\subset \tim_*.$
 
 \ss \item[(c)] Under the hypotheses of part {\rm (b)}, if 
  $\mu_X$  is supported in an open ball $B$ of 
 radius less than a certain number $\rcx(M,g)$ {\em (specified in Definition \ref{def:rinj_rcx})},
then almost surely the Fr\'{e}chet mean set of $\{X_1,\dots,X_N\}$
is a single point and lies in  $M_*\cap B$. Correspondingly, under the same hypotheses,
almost surely 
 the equivariant Fr\'{e}chet means of $\{X_1,\dots,X_N\}$ are unique
up to the action of $G$ and lie in $\tim_*\cap\quo^{-1}(B)$. 
\end{itemize}
 \end{thm}

Our strategy for proving
Theorem \ref{thm1.1} is to deduce it from a purely geometric  (and quite general) result that
constitutes the core mathematics in this paper:

\begin{thm}\label{thm:FMmeaszero}
If $A\subset M$ has volume zero in $M$ , then
the set
\be
\{Q\in M^N : \FM(Q) \intersect A\neq\emptyset\}
\ee
has volume zero in $M^N$. In particular, this holds if $A$ is a countable union of submanifolds
 of positive codimension.
\end{thm}

\ss 
Since Theorem \ref{thm1.1}  is derived as a consequence of  Theorem \ref{thm:FMmeaszero},
much of this paper is, correspondingly, purely geometric.  Our approach is to turn the question
answered by Theorem \ref{thm:FMmeaszero} into a {\em transversality} problem by using
``pre-barycenters''.  These are pre-images of barycenters in the bundle $\tmn$, the $N$-fold fiberwise product of $TM$ with itself.  (We use the term {\em barycenter}  in the  ``balancing point'' sense in this paper.)
The intuition behind Theorem \ref{thm:FMmeaszero} is that if we perturb a sample $Q_0$ to
a sample $Q$ that moves over an
open neighborhood of $Q_0$---or even if we hold all but one of the sample points fixed, and allow the remaining
point to vary over an open neighborhood---the barycenters should {\em move}, varying over an open set themselves.
In Sections \ref{sec:bary-prebary}--\ref{sec:barymanifold} we define and discuss barycenters and pre-barycenters; show that the set $\pb_N$ of {\em pre-barycentered configurations} is a submanifold of the bundle $\tmn$ (part of Proposition \ref{more_submersions}(a)); prove a forerunner of Theorem \ref{thm:FMmeaszero}, a genericity property of barycenters rather than Fr\'{e}chet means (Theorem \ref{thm:measzero.7}); and establish some properties of the set of $\B_N$ of {\em barycentered configurations}
(the set $\B_N := \{(Q,p): \mbox{$p$ is a barycenter of $Q$}\}\subset M^N\times M\}$). 

Among the latter properties is that the set of {\em sufficiently concentrated} 
barycentered configurations is a submanifold of $M^N\times M$ (part of Corollary \ref{cor:conFMbary}), a property that we do not know or expect of the entire set $\pb_N$ in general. In Section \ref{sec:FMs_are_barys}, 
 building on work in \cite{le-barden2014},  we show that all Fr\'{e}chet means of finite configurations are barycenters (Corollary \ref{cor:fmisbary}).  This fact, combined with Theorem \ref{thm:measzero.7}, yields Theorem \ref{thm:FMmeaszero}.  

In Section \ref{sec:iid1} we begin applying our geometric results to i.i.d. random variables. 
 Theorem \ref{thm1.1}
follows quickly from Theorem \ref{thm:FMmeaszero}.  
As mentioned earlier, we apply this together with \cite{afsari}'s Theorem 2.1 to obtain, additionally, the almost-sure uniqueness of Fr\'{e}chet sample means (under appropriate hypotheses). In Section \ref{sec:iid1} we also discuss an example that illustrates the generality of Theorem \ref{thm1.1}, and also give some examples of interest in which 
Theorem \ref{thm1.1} does {\em not} apply. 

Section \ref{sec:eqvt_means_general} is pure geometry again:  we examine finite-degree principal Riemannian 
covers $(\tim, \tg)\stackrel{\quo}{\longrightarrow}(M,g)$ and equivariant means.   For a finite configuration $Q$, we show that taking equivariant Fr\'{e}chet means in $\tim$ is ``equivalent''
to taking
Fr\'{e}chet means in $M$, in the sense that $\EFM(Q)=\quo^{-1}(\FM(Q))$ (Proposition \ref{prop:uniqueuptoG}(a)).   We also establish relations between the numbers
$\rcx(\tim,\tg)$ and  $\rcx(M,g)$.
We are then set up to 
address Fr\'{e}chet sample means in the equivariant setting in Section \ref{sec:iid_eqvt_means}, 
where we prove Theorem \ref{thm1.2} (as Corollary \ref{cor:measzero.5}).

In Section \ref{sec:PSRmeans} we discuss an application of our equivariant-Frechet-means results to ``partial scaling-rotation means'', a type of mean for symmetric positive-definite matrices that the authors introduced in \cite{rjgs}.

Some remarks concerning the terminology ``volume zero'', ''generic'', ``absolutely continuous'', 
and ``$M$-valued random variable'' are worthwhile here.  
The notion of {\em measure-zero subset of a manifold $\M$} (see e.g. \cite[Ch. 6]{lee_mfds} or \cite[Ch. 3]{hirsch})
is  independent of any measures that one is putting on $\M$, but coincides with ``set of ($n$-dimensional) Lebesgue measure zero''  when $\M=\bfr^n$. 
%{\color{red} [DG note: ``$\M$'' instead of $M$ is not a typo.  Here I want
%to talk about a general manifold that need not be the one on p. 1.]}
To avoid any
confusion with probability measures under consideration, we will refer to 
a measure-zero subset of a manifold
$\M$ (in the sense above) 
 as having {\em volume zero} in $\M$; we call the complement of such a subset {\em generic} (in $\M$). 
Motivation for the ``volume zero'' terminology is that
every Riemannian metric $g$ on $\M$
determines a {\em volume measure} whose null sets are precisely the measure-zero sets in the sense above. 
However, the above notions of ``volume zero'' and ``genericity''---{\em residuality} in the terminology of \cite[Ch. 3]{hirsch}---are at heart more topological than measure-theoretic, making sense even without specifying a measure on $M$.
Every generic subset of $\M$ is dense, and every positive-codimension submanifold of $\M$
of  has volume zero in $\M$.

We call a measure $\mu$ on a manifold $\M$  {\em absolutely continuous} if $\mu(A)=0$
whenever \linebreak $A\subset \M$ has volume zero in $\M$. Every measure on $\M$ given by
a continuous density function is absolutely continuous. 

Every manifold $\M$ has two ``natural'' sigma-algebras of (potentially) measurable sets:  the sigma-algebra ${\mathcal B}_\M$ of Borel sets, and its completion ${\mathcal L}_\M$, the ``Lebesgue sigma-algebra''---an analog of the sigma-algebra of Lebesgue measurable subsets of $\bfr^n$---which can be obtained by adding to ${\mathcal B}_\M$ 
the collection of all volume-zero subsets of $\M$. 
To have a clear meaning of  ``measurable map''
(a map for which the inverse image of every measurable set is measurable),
with domain or codomain $\M$,   we must specify which of these sigma-algebras on $\M$ we are using.   
It is desirable to make the same specification consistently, across all manifolds, to ensure that the
composition of measurable maps is measurable.  
Since volume-zero sets occur throughout this paper, we choose the Lebesgue sigma-algebra
on every manifold. (This choice differs from the customary definition of ``measurable function $f:\bfr\to\bfr$'', 
in which $f$ is called measurable if the inverse image of every {\em Borel} set is {\em Lebesgue} measurable.)

In this paper we will take the domain of every random variable that appears to be a (fixed but arbitrary) complete
probability space
$(\Omega,\mathcal{A},\mathcal{P})$.
An $\M$-valued random variable is then an $(\mathcal{A},{\mathcal L}_\M)$-measurable map 
$X:\Omega\to \M.$ We
call  $X$ {\em absolutely continuous} if its induced measure $\mu_X$ on $\M$ is absolutely continuous.  For such $X$, it actually makes no difference whether we choose ${\mathcal L}_\M$ or ${\mathcal B}_\M$ as our sigma-algebra on $\M$: if $X:\Omega\to \M$ is an $(\mathcal{A},{\mathcal B}_\M)$-measurable map 
whose induced measure $\mu_X$ on $\M$ is absolutely continuous, then $X$ is also 
$(\mathcal{A},{\mathcal L}_\M)$-measurable.

Throughout this paper, $n=\dim(M),$ and $N$ is an arbitrary positive integer.
For
$Q=(q_1,\dots,q_N)\in M^N$ we write
$\supp(Q)$ for the (unordered) multi-set $\{q_1,\dots,q_N\}$; when the points are distinct
this coincides with the support of the measure $\mu_Q$.
For $p\in M$, the zero element of $T_pM$ is sometimes denoted $0_p$ for clarity's sake, and
sometimes simply denoted $0$.
As is customary, ``$\exp$'' denotes the exponential map $TM\to M$, while
``$\exp_p$'' denotes the restriction of $\exp$ to $T_pM$.
In any metric space that is understood from context, $B_r(q)$ and $\bar{B}_r(q)$
denote, respectively, the open and closed balls of radius $r$ centered at $q$; sometimes
this notation may be augmented to include some reference to the metric space of
interest.

 To avoid technical distractions, throughout we take ``smooth'' to mean $C^\infty$, and ``(sub)manifold'' 
to mean ``smooth (sub)manifold''.  These restrictions of ``smooth'' and  ``(sub)manifold'' can be relaxed; 
see Remark \ref{rem:findiff}.

\setcounter{equation}{0}
\section{Barycenters and pre-barycenters}
\label{sec:bary-prebary}

In this paper, we take {\em barycenter} of a configuration $Q=(q_1,\dots, q_N)\in M^N$
to mean a ``balancing point'' (suitably defined) for the configuration.
 The simplest example
of such a barycenter is a point $p$ for which no $q_i$  lies in the cut-locus of $p$
(thereby allowing a reasonable definition of ``$(\exp_p)^{-1}(q_i)$'') and for which
$\sum_{i=1}^N (\exp_p)^{-1}(q_i) =0$.
This example is what we will soon be calling a {\em short barycenter}.
To give a precise and more general definition of ``barycenter'' and related notions,
we recall several facts pertaining to the exponential map and cut-loci, and fix
some notation.   All these facts can be found in standard books on Riemannian geometry, e.g. \cite{chavel,sakai}, except for those facts for which we give specific references.

Let $S(TM)$ denote the unit tangent bundle of $(M,g)$, and for $p\in M$ let $S(T_pM)=S(TM)\intersect T_pM$.  
For $v\in TM$ let $\g_v$ be the geodesic $t\mapsto \exp(tv)$. The {\em cut-time} function
$\tau:S(TM)\to [0,\infty]$ is defined by $\tau(v)=\sup\{t\geq 0: \g_v|_{[0,t]}\ \mbox{is minimal}\}$.
Recall that the {\em tangent cut locus} $\CLtan(p)$ and {\em cut locus} $\CL(p)$ of $p\in M$
are defined by $\CLtan(p)=\{\tau(v)v: v\in S(T_pM), \tau(v)<\infty\}$
and $\CL(p)=\exp_p(\CLtan(p))$. Following \cite{bishop}, we call
a cut-point $q\in \CL(p)$ {\em ordinary} or {\em singular} accordingly as there are at least two minimal geodesics
from $p$ to $q$, or only one, and classify $v\in \CLtan(p)$ as an ordinary or singular accordingly as
$\exp_p(v)$ is an ordinary or singular cut-point of $p$.
We will write $\clord(p)$ and $\clsing(p)$ for the set of ordinary and singular cut-points of $p$, respectively, and
$\clsingtan(p)$ for the set of singular points in the tangent cut-locus of $p$.
Note that any of these sets can be empty, and that the relations ``$q\in \clord(p)$'' and ``$q\in\clsing(p)$'' are symmetric in $p$ and $q$.

Still for arbitrary $p\in M$, define 
\bestar
\hspace{-1.25in}
\Dtan(p)&=&\{tv: v\in S(T_pM), t\in [0,\tau(v))\},\\
\tDtan(p)&=&\Dtan(p)\,\union\,\clsingtan(p),\\
\D(p)&=&\exp_p(\Dtan(p)), \\ 
\mbox{and}\hspace{.75in} \tD(p)&=& \exp_p(\tDtan(p)) \ =\ \D(p)\,\union\,\clsing(p).
\hspace{.75in}
\eestar

\ssn
(Thus $\tD(p)$ is precisely the set of points  $q\in M$ for which
there is a unique minimal geodesic from $p$ to $q$.)
The set $\Dtan(p)$ is open, and the map $\exp_p$ is injective on $\tDtan(p)$ and
nonsingular on $\Dtan(p)$.
The set $\tDtan(p)$ is the largest subset $U\subset T_pM$, star-shaped
with respect to $0_p$, such that $\exp_p|_U$ is injective; $\Dtan(p)$ is the largest {\em open}
subset $U\subset T_pM$, star-shaped
with respect to $0_p$, such that $\exp_p|_U$ is a diffeomorphism onto its image. In particular,
$\D(p)$ is open, $\exp_p|_{\Dtan(p)}:\Dtan(p)\to \D(p)$ is a diffeomorphism, and
$M=\D(p)\,\disjoint\,\CL(p)$.
We define the notation  $(\exp_p)^{-1}$ to mean the inverse
of the
bijective map $\exp_p|_{\tDtan(p)}: \tDtan(p)\to\tD(p)$. Note that $(\exp_p)^{-1}$ may
not be continuous at points of $\clsing(p)$, but the restriction of
$(\exp_p)^{-1}$ to the open set $\D(p)$ is a diffeomorphism.  The open set $\D(p)$ is
dense in $M$, and its complement $\CL(p)$ has measure zero.
The distance function $r_p$ is never $C^1$ at an ordinary cut-point of $p$, and
never $C^2$ at any cut-point of $p$ \cite{bishop}, but $r_p^2$ is smooth 
on $\D(p)$.

Our definition of {\em barycenter} is now as follows.

\begin{defn}\label{def:bary} \rm Let $Q=(q_1,\dots,q_N)\in M^N$ and let $p\in M$.  We call $p$ a {\em barycenter
of $Q$} if there exist $v_1,\dots, v_N\in T_pM$ such that $q_i=\exp_p(v_i), 1\leq i\leq N$, and
$\sum_i v_i=0$.  We call $p$ a {\em short} barycenter of $Q$ if there exist such $v_1,\dots,v_N
\in
\D^{\tan}(p)$, and an {\em almost-short} barycenter of $Q$ if there exist such $v_1,\dots,v_N
\in
\tD^{\tan}(p)$.
We will write $\B(Q)$ for the set of barycenters of $Q$, and $\B^0(Q)$ (respectively, $\tB^0(Q)$)
for the set of short (resp., almost-short) barycenters of $Q$.  We define
\ben
\B_N
=\B_N(M)=\{(Q,p)\in M^N\times M \ : \  p\in
\B(Q)\},
\een
the set of {\em barycentered configurations of $N$ points in $M$}. Similarly, we define
$\B_N^0=\{(Q,p)\in M^N\times M \ : \  p\in \B^0(Q)\}.$
Note that $\B(Q), \B^0(Q)$, and $\tB^0(Q)$ are subsets of $M$, while
$\B_N$ and $\B_N^0$ are  subsets of $M^N\times M$.
\end{defn}

Also important will be the notion of {\em pre-barycenter}.  To define this,
let $\tmn$ denote the fiber bundle over $M$ whose fiber over $p\in M$ is the $N$-fold
Cartesian product $(T_pM)^N$.  (The bundle-structure is the same as that of the
$N$-fold Whitney sum of $TM$ with itself, but we are regarding the fiber simply as a Cartesian product rather than a direct sum.)  We will write $\pi:TM\to M$ and $\pi':\tmn\to M$ for the projection-maps
of these bundles. Thus, elements of $\tmn$ may be written as $V=(v_1,\dots, v_N)
\in TM\times \dots\times TM$ with $\pi(v_1)=\pi(v_2)=\dots=\pi(v_N)$, and
we have $\pi'(V)=\pi(v_i), 1\leq i\leq N$. Note that $\dim(\tmn)=(N+1)n=\dim(M^N\times M)$.

\begin{defn}\label{def:prebary}\rm Define
\ben
\pb_N=\pb_N(M)=\{\ (v_1,\dots,v_N)\in\tmn : \sum_i v_i=0\ \}.
\een
We call elements of
$\pb_N$ {\em pre-barycenters} (of configurations of $N$ points in $M$).
We define the map $\exp^N:\tmn\to M^N$ by $\exp^N(v_1,\dots,v_N)=(\exp(v_1),\dots,\exp(v_N))$.
 If $V= (v_1,\dots,v_N)\in \pb_N$ and
$Q=\exp^N(V)$, we say that $V$ is a pre-barycenter of $Q$. For $Q\in M^N$ we
write $\pb(Q)$ for the set of pre-barycenters of $Q$.
\end{defn}

From Definitions \ref{def:bary} and \ref{def:prebary} it is clear that for all $Q\in M^N$,
\be\label{baryprebary}
\B(Q) = \pi'(\pb(Q)).
\ee

\begin{remark} \rm
We are using {\em ordered} $N$-tuples $Q=(q_1,\dots,q_N)$ to take advantage of the manifold
structures
of $M^N$ and $\tmn$.  Clearly the set
$\B(Q)$
depends only
on the measure $\mu_Q$, not on how the points $q_i$ are ordered.
\end{remark}

\begin{remark}\rm Our definitions of {\em barycenter} and {\em pre-barycenter} of a configuration
are closely related to the
definition of {\em barycentre of a (Borel) probability measure} in \cite{corcuera}, which
for purposes of distinction with our definition
we will call a CK-barycenter.  Given a measure $\mu$ on $M$ and a point $p\in M$, we call
a measure $\tilde{\mu}$ on $T_pM$ a {\em lift} of $\mu$ if $\mu$ is the push-forward,
under $\exp_p$, of  $\tilde{\mu}$. A point $p\in M$ is a CK-barycenter of a probability measure $\mu$ if
$\mu$ lifts to some probability measure $\tilde{\mu}$ on $T_pM$ having zero  ``vector mean'',
i.e. satisfying $\int_{T_pM} v\,d\tilde{\mu}(v)=0$.
For $V=(v_1,\dots,v_N)\in \tmn$, $p=\pi'(V)$, and $Q=\exp^N(V)$, the probability measure
$\mu_V:=\frac{1}{N}\sum_i \d_{v_i}$ on $T_pM$ is a lift of $\mu_Q$.  Thus if $V$ is a pre-barycenter
of $Q$, or $p$ is a barycenter of $Q$, then $p$ is a CK-barycenter of $Q$.  Our pre-barycenters of $Q$ are a proxy for
the probability measures of the form $\mu_V$ on $T_pM$ to which $\mu_Q$
may lift.   However, our definition of ``barycenter of $Q$'' may be more restrictive than
``CK-barycenter of $\mu_Q$'',  since for $p$ to be a barycenter of $Q$ we require that $\mu_Q$  lift not just
to an arbitrary vector-mean-zero probability measure on $T_pM$, but to one
of the form $\mu_V$ above. (A general lift of $\mu_Q$ need not even be discrete; if discrete,
it may have more than $N$ points, and the weights of certain points need not be equal.)
\end{remark}

Taking the liberty of writing elements of $M^N\times M$ as elements of $M^{N+1}$,
we define subsets $\tM_1, \M_1$ of $M^N\times M$ by
\bearray\label{defm4}
\tM_1&=&\{(q_1,\dots,q_N,p)\in M^N\times M \ : \ 
q_i\in \tD(p), \ 1\leq i\leq N\},\\
\M_1&=&\{(q_1,\dots,q_N,p)\in M^N\times M \ : \ 
q_i\in \D(p), \ 1\leq i\leq N\},
\eearray
and define a map $F:\tM_1\to TM$ by
\bearray\nonumber
F(q_1,\dots, q_N,p) &=&\frac{1}{N}\sum_{i=1}^N (\exp_p)^{-1}(q_i)\ \ \  \in T_pM.
\label{defFQp}
\eearray
The set $\M_1$ is open and dense in $M^N\times M$, and the restricted map $F|_{\M_1}$ is smooth.

\begin{notation}\rm {\ }
\begin{enumerate}
\item[(a)]
Let $\pi_N : M^N\times M\to M^N, \plast : M^N\times M\to M$ denote
projection onto the first and second factors of $M^N\times M$. We use the same notation for the restrictions of $\pi_N, \plast$ to  $\M_1$ and $\tM_1$.  Observe that $\plast(\M_1)=M$.

\item[(b)] For $Q=(q_1,\dots, q_N)\in M^N$, define the following:

\begin{enumerate}
\item[(i)] $f_Q=\frac{1}{N}\sum_{i=1}^N r_{q_i}^2: M\to\bfr$.

\item[(ii)]
$M_Q=\{p\in M: (Q,p)\in \M_1\}=\{p\in M\ : \ 
\supp(Q)\subset \D(p)\}= 
\plast(\pi_N^{-1}({Q})\cap \M_1)$.  Since the condition $\supp(Q)\subset \D(p)$
is equivalent to $p\in \Intersect_{i=1}^N \D(q_i)$, and since $\D(q_i)$ is open and dense in $M$
for each $i$, the set $M_Q$ is open and dense as well.

\item[(iii)]
$\tM_Q=\{p\in M: (Q,p)\in \tM_1\}=\{p\in M\ : \ 
\supp(Q)\subset \tD(p)\}= 
\plast(\pi_N^{-1}({Q})\cap \tM_1)$.

\item[(iv)] Define a vector field $Y_Q$ on $M_Q$ by $Y_Q(p)=F(Q,p)$.

\end{enumerate}

\end{enumerate}

\end{notation}

It is well known that $f_Q$ is smooth on $M_Q$ and that the vector field $Y_Q$ is the negative gradient of $f_Q|_{M_Q}$
(see e.g. \cite{karcher}).  Thus the following are equivalent: (i) $p$ is a short barycenter
of $Q$; (ii) $p$ is a critical point
of $f_Q$ in $M_Q$; and (iii) $Y_Q(p)=0$.
  In particular,  $\B^0(Q)$ contains all  Fr\'{e}chet means
of $Q$ lying in $M_Q$.

We will see later that $\B(Q)$ contains the entire set $\FM(Q)$; in fact
 $\FM(Q)\subset
\tB^0(Q)$. This nontrivial fact, which follows from results in \cite{le-barden2014},
will essentially allow us to deduce Theorem \ref{thm1.1} from a stronger conclusion in which
$\FM(\{X_1,\dots,X_N\})$ in \eqref{eq:thm1.1} is replaced by $\B(\{X_1,\dots,X_N\})$.
For efficiency of presentation we structure the argument slightly differently, restating
and proving
Theorem \ref{thm1.1} as Corollary \ref{cor:measzero.3} (a).

\setcounter{equation}{0}
\section{Genericity of barycenters}
\label{sec:bary-prebary2}

At the heart of Theorem \ref{thm1.1} is the question of how the
set $\FM(Q)$
varies as we vary $Q$.  Rather than approaching this
directly, our strategy will be to analyze
perturbations of
(pre-)barycenters and to show that from this analysis, the desired facts about
Fr\'{e}chet means can be deduced.  Perturbation analysis is easier for barycenters, short
barycenters, and pre-barycenters than (directly) for
Fr\'{e}chet means since barycenters (and short and pre-barycenters) are defined as solutions of {\em equations} on open subsets of manifolds, putting us in the realm in which we can hope to
apply the Implicit Function Theorem in
some guise (e.g. the Regular Value Theorem or transversality). This is why we have
defined the sets $\B_N, \B_N^0,$ and $\pb_N$---parametrized families of
(pre-)barycenters, parametrized by the $N$-tuple $Q$.  As we shall see shortly, $\pb_N$
is a manifold, and the family of barycenter-sets $\B(Q)$ enjoys a genericity property.

\begin{figure}[t]
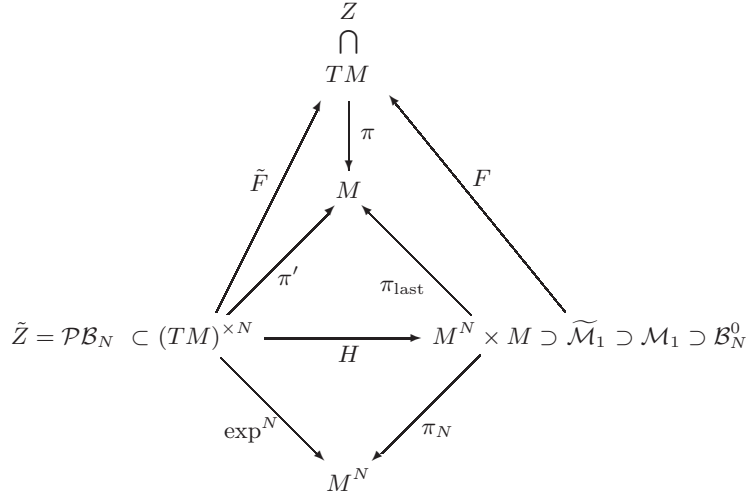

\begin{tabular}{c}

\begin{diagram}[nohug]
 & & &
{
\begin{array}{c}
Z\\
\bigcap
\\ \phantom{\int_a^{b^c}}\mbox{\hspace{-.2in}}TM
\end{array}
}
& & & &\\
& & \ruTo(2,4)<\tF &
\dTo>{\mbox{\hspace{-.07in}}\begin{array}{l}~\\ \,\pi\end{array}}&\luTo(3,4)>F &  &&\\
& & & M & &  &&\\
& &\ruTo>{\pi'}& &\luTo<\plast & &&\\
\tz=\pb_N \ \subset &\ \tmn && \rTo_{{H}} & & {M^N\times M} \supset &\  \tM_1\
&\supset \M_1 \supset \B_N^0
\\
& & \rdTo<{\exp^N} & &\ldTo>{\pi_N} & &&\\
& & & M^N & &  &&\\
\end{diagram}

\end{tabular}

\caption{Commutative diagram depicting relations among various spaces and maps.}
\label{comdiag2}
\end{figure}

To this end, writing elements of $\tmn$ as $V=(v_1,\dots,v_N)$, we define maps $\tF:\tmn\to TM$ and
${H}:\tmn\to M^N\times M$ by
\bestar
\tF(V)&=&\sum_i v_i\ ,\\
{H}(V)&=&(\exp^N(V), \pi'(V)).
\eestar
Henceforth let $Z$ denote the zero-section of $TM$ (as a submanifold of $TM$).  Observe
that
\be\label{finvz-2}
\B_N^0= (F|_{\M_1})^{-1}(Z)\ \ \ \mbox{and}\ \ \
\pb_N=\tF^{-1}(Z)=: \tz.
\ee
(We use the notation $\tz$ for $\pb_N$ to simplify other notation below.)
The commutative diagram in Figure \ref{comdiag2} depicts relations among various
spaces and maps we have defined.  Note that if we define
\be\label{defm4'}
\M_1'=\{V=(v_1,\dots,v_N)\in \tmn \ : \ 
v_i\in \Dtan(\pi'(V)), \ 1\leq i\leq N\},
\ee
then $\tF|_{\M_1'\cap \tz}$ is the inverse of the natural injection
\bearray\nonumber
\B_N^0 &\to& \tz, \\
(q_1,\dots, q_N,p) & \mapsto & ((\exp_p)^{-1}(q_1), \dots, (\exp_p)^{-1}(q_N)).
\label{includebary}
\eearray

\bs
In the next proposition and beyond, recall that $n=\dim(M)$.

\begin{prop} \label{more_submersions}
(a) The map $\tF$ is a submersion, and $\tz$---the set of pre-barycentered configurations of $N$ points 
in $M$---is a submanifold of $\tmn$ of dimension $Nn$.

(b)
The map $\pi'|_{\tz}:\tz\to M$ is a submersion.
\end{prop}

The horizontal/vertical splitting of the tangent spaces of $TM$ and $\tmn$, induced by the Riemannian
structure, enters the proof of Proposition \ref{more_submersions}, so we review these splittings before
proceeding with the proof.

For each $v\in TM$, the {\em vertical space} $\Ver_v(TM)\subset T_v(TM)$ is defined to be
$\ker(\pi_{*v})$.  Letting $p=\pi(v)$, there is a canonical isomorphism $\iota_v: T_pM\to \Ver_v(TM)$
defined by  $\iota_v(w)=\frac{d}{dt}(v+tw)|_{t=0}$.  The Levi-Civita connection $\na$
determines a complement of $\Ver_v(TM) \subset T_v(TM)$, the {\em horizontal space} $\Hor_v(TM)$,
also isomorphic to $T_pM$; the map $\pi_{*v}|_{\Hor_v(TM)}:\Hor_v(TM)\to T_pM$ is an
isomorphism. Given a curve
based at $p$---a smooth map
$\g:I\to M$ where $I$ is an interval containing $0$ and where $\g(0)=p$---the {\em horizontal lift}
of $\g$ starting at $v$ is the curve $\tgam:I\to TM$ defined by $\tgam(t)=\mathcal{P}_{\g(0)\to
\g(t)} v$, the parallel transport of $v$ along $\g$ from time $0$ to time $t$.  Given $X\in T_pM$,
the {\em horizontal lift} of $X$ to $TM$ at $v$ is the vector $\tilde{X}_v=\tgam'(0)$, where $\g$ is any curve with $\g'(0)=X$.
The space $\Hor_v(TM)$ is the set of all horizontal lifts of elements of $T_pM$, and the map $X\mapsto\tilde{X}_v$ is the inverse of the isomorphism $\pi_{*v}|_{\Hor_v(TM)}$.

The Levi-Civita connection similarly induces a horizontal/vertical splitting of the tangent spaces of $\tmn$. Let $V=(v_1,\dots,v_N)\in \tmn$ and let $p=\pi'(V)$.  Viewing $\tmn$ as the restriction to the
diagonal of $N$-fold Cartesian product of $TM$ with itself (a bundle over $M^N$), the tangent space $T_V(\tmn)$ may be canonically identified with a subspace of
$\Plus_{i=1}^N T_{v_i}(TM)$, namely $\{(w_1,\dots, w_N)\in \Plus_{i=1}^N T_{v_i}(TM)
\ : \  \pi_{*v_1}w_1 = \pi_{*v_2}w_2=\dots=\pi_{*v_N}w_N\}$.  The vertical subspace of
$T_V(\tmn)$ is

\bestar
\Ver_V(\tmn):=\ker(\pi'_{*V})
&=&
\Plus_{i=1}^N\Ver_{v_i}(T_pM)\\
&=&
\{(\iota_{v_1}(w_1), \dots, \iota_{v_N}(w_N)) \ : \  w_i\in T_pM\}\\
&\caniso& \Plus_{i=1}^N T_pM.
\eestar
\noi For $X\in T_pM$, we define
the horizontal lift of $X$ to $\tmn$ at $V$ to be the vector
\be\label{txv}
\tilde{X}_V=
(\tilde{X}_{v_1}, \dots, \tilde{X}_{v_N})
\in T_V(\tmn)\subset \Plus_{i=1}^N T_{v_i}(TM),
\ee
 where
$\tilde{X}_{v_i}$ is the horizontal lift of $X$ to $TM$ at $v_i$, $1\leq i\leq N$.
Defining $\Hor_V(\tmn)$ to be the set of all horizontal lifts of elements of $T_pM$ to $\tmn$,
it is easily seen that $\Hor_V(\tmn)$ is a vector subspace of $T_V(\tmn)$ and is a complement
to $\Ver_V(\tmn)$.

\bsn
{\bf Proof of Proposition \ref{more_submersions}.}
Let $V=(v_1,\dots,v_N)\in \tmn$ and let $p=\pi'(V)$.  Decomposing $T_V(\tmn)$ as
$\Hor_V(\tmn)\plus \Ver_V(\tmn)$, let us write elements of $T_V(\tmn)$ in the form
$(\tilde{X}_V, (\iota_{v_1}(w_1), \dots, \iota_{v_N}(w_N)))$, where $\tilde{X}_V$ is the horizontal
lift at $V$ of $X\in T_pM$ and $w_1,\dots,w_N\in T_pM$.  Then one easily computes
\be\label{tf_*V}
\tF_{*V}(\tilde{X}_V, (\iota_{v_1}(w_1), \dots, \iota_{v_N}(w_N)))
=(\tilde{X}_{\tF(V)}, \iota_{\tF(V)}(\sum_i w_i))
\ee
where $\tilde{X}_{\tF(V)}\in \Hor_{\tF(V)}(TM)$ is the horizontal lift of $X$ to $TM$ at $\tF(V)$.
Clearly the image of $\tF_{*V}$ is all of $\Hor_{\tF(V)}(TM)\plus \Ver_{\tF(V)}(TM) = T_{\tF(V)}(TM)$.
Hence $\tF$ is a submersion, and $\tz=\tF^{-1}(Z)$ is a submanifold of $\tmn$ of codimension
equal to the codimension of $Z$ in $TM$, namely $n$.  Hence
$\dim(\tz)=\dim(\tmn)-n =Nn$.

 Let $V=(v_1,\dots,v_N)\in \tF^{-1}(Z)$ and $p=\pi'(V)$; thus $\tF(V)=0_p\in T_pM$. Recall that $\Hor_{0_p}(TM)=T_{0_p}Z$.
Hence using \eqref{tf_*V} we find
\bearray\nonumber
\lefteqn{T_V(\tz)=(\tF_{*V})^{-1}(T_{\tF(V)}Z)=(\tF_{*V})^{-1}(\Hor_{0_p}(TM))}\\
&=& 
\{(\tilde{X}_{\tF(V)}, (\iota_{v_1}(w_1), \dots, \iota_{v_N}(w_N)))\in T_V(\tmn)\ : \  \sum_i w_i =0_p\}.
\label{TVGinvZ}
\eearray
Note that $\Hor_{\tF(V)}(T_V(\tmn))\subset T_V\tz$ and that $\pi'_{*V}(\Hor_V(\tmn))
=T_pM$; thus
$\pi'_{*V}|_{T_V(\tz)}$ is surjective.  Hence $\pi'|_{\tz}$ is a submersion.
\qedns

Knowing that $\pi'_{\tz}$ is a submersion (Proposition \ref{more_submersions}(b)) will be valuable to us because of the following lemma.  This lemma is almost certainly known, but the authors have found no reference for it. For a proof, see
Appendix \ref{sec:appendixA}.

\begin{lemma}\label{submersions_and_meas_zero}
Under a submersion, the inverse image of a set of volume zero has
volume zero.
\end{lemma}

\begin{thm}\label{thm:measzero.7}
If $A\subset M$ has volume zero in $M,$ then
the set
\be\label{FMinverse(A),2}
\{Q\in M^N : \B(Q) \intersect A\neq\emptyset\}
\ee
has volume zero in $M^N$.
\end{thm}

\pf Let $A$ be a volume-zero subset of $M$.  By Proposition \ref{more_submersions}
and Lemma \ref{submersions_and_meas_zero}, $(\pi'|_{\tz})^{-1}(A)$ has volume zero in
the $Nn$-dimensional manifold $\tz$.  Every smooth map between manifolds  of the same
dimension carries sets of volume zero to sets of volume zero \cite{lee_mfds}.
In particular, this is true of the map $\pi_N\circ {H}|_{\tz}:\tz\to M^N$.
Thus
$\pi_N({H}((\pi'|_{\tz})^{-1}(A)))$ has volume zero in $M^N$.
Using Corollary \ref{cor:fmisbary}, observation \eqref{baryprebary}, and the fact that
$\plast\circ{H}=\pi'$,
\bearray
\nonumber
\{Q\in M^N \ : \  \B(Q) \intersect A\neq\emptyset\}
\nonumber
&=&
\{Q\in M^N \ : \  \exists p\in A\ \mbox{such that}\ p\in \B(Q)\}\\
\nonumber
&\subset&
\{Q\in M^N \ : \  \exists p\in A\ \mbox{such that}\ p\in \plast(\pb(Q))\}\\
\nonumber
&=&
\{Q\in M^N \ : \  \exists p\in A\ \mbox{such that}\ (Q,p)\in {H}(\tz)\}\\
\nonumber
&=&
\pi_N\left({H}(\tz)\intersect \plast^{-1}(A)\right)\\
&=&
\pi_N\left({H}\left((\pi'|_{\tz})(A)\right)\right).
\label{volzeroset}
\eearray
Hence the set \eqref{FMinverse(A),2} has volume zero in $M^N$.\qedns

\setcounter{equation}{0}
\section{Some properties of the set of barycentered configurations}
\label{sec:barymanifold}

Proposition \ref{more_submersions}(a) shows that $\pb_N$, the set of {\em pre}-barycentered configurations, is a submanifold
of $\tmn$.  It appears unlikely that $\B_N$, the set of barycentered configurations, is
a submanifold of $M^N\times M$.  However, we shall see that the subset $\B_N^0$
{\em is} a submanifold of $M^N\times M$, and captures  the barycenters of almost all
configurations in the
sense that
$\pi_N(\B_N\setminus\B_N^0)=M^N\setminus\pi_N(\B_N^0)$ has volume zero
in $M^N$.  (We have
$\pi_N(\B_N)=M^N$ since $\FM(Q)\neq\emptyset$ for every $Q\in M^N$.)

In proving these facts,
two types of derivative of the vector fields $Y_Q$ will arise.  One
is the covariant derivative $\na Y_Q$, whose value at $p\in M$ is an  endomorphism of $T_pM$,
namely the map
$v\mapsto \na_v Y_Q$.  The other is the derivative of the map-of-manifolds $Y_Q:M\to TM$, whose
value at $p$ is the linear map $Y_{Q*p}:T_pM\to T_{Y_Q(p)}(TM)$.  These derivatives are related
as follows: for $v\in T_pM$,
\be\label{YQ*p}
Y_{Q*p}(v)=\tilde{v}_{_{Y_Q(p)}} + \iota_{_{Y_Q(p)}}(\na_v Y_Q);
\ee
i.e. $\na_v Y_Q$ is essentially the vertical part of $Y_{Q*p}(v)$ \cite{groisser2004, karcher}.

\begin{lemma}\label{lem:F_*}
 Let $(Q,p)=(q_1,\dots,q_N,p)\in \M_1$.  Writing elements of $T_{(Q,p)}(M^N\times M)$
in the form $(w_1,\dots, w_N,v)$, where $w_i\in T_{q_i}M$ and $v\in T_pM$,
the derivative of $F$ at $(Q,p)$ is given by
\bearray\nonumber
\lefteqn{F_{*(Q,p)}(w_1,\dots, w_N,v)=}
\\
&&
\tilde{v}_{_{Y_Q(p)}}+
\iota_{_{Y_Q(p)}}\left\{ \na_v Y_Q +
\frac{1}{N}\sum_{i=1}^N \iota_{z_i}^{-1}\left[((\exp_p)^{-1})_{*q_i}(w_i)
\right]
\right\},
\label{F*qp.0}
\eearray
\noi where $z_i=(\exp_p)^{-1}(q_i)$.
\end{lemma}

\pf Let $(w_1,\dots, w_N,v)\in T_{(Q,p)}\M_1=T_{(Q,p)}(M^N\times M)$. We have
\bearray
\nonumber
F_{*(Q,p)}(w_1,\dots, w_N,v)&=&
\sum_{i=1}^N F_{*(Q,p)}(0_{q_1}, \dots, 0_{q_{i-1}}, w_i,0_{q_{i+1}},\dots,0_{q_N},0_p)\\
&&+F_{*(Q,p)}(0_{q_1},\dots, 0_{q_N},v).
\label{F*Qp}
\eearray
For $1\leq i\leq N$, the vector $((\exp_p)^{-1})_{*q_i}(w_i)$ is vertical:  it is the initial
velocity of a curve of the form $t\mapsto
(\exp_p)^{-1}(\g_i(t))$ lying entirely in $T_p M\subset TM$ (where $\g_i$ is a
curve in $M$ defined on some interval $(-\e,\e)$ and having $\g_i'(0)=w_i$).
The $i^{\rm th}$ summand on the first line of \eqref{F*Qp} is simply the image
of $((\exp_p)^{-1})_{*q_i}(w_i)$ under the isomorphism $\iota_{Y_Q(p)}\circ \iota_{z_i}^{-1}$;
the isomorphisms $\iota_{z_i}^{-1}$ allow us to identify the vertical vectors
 $((\exp_p)^{-1})_{*q_i}(w_i)\in T_{z_i}M$ with elements of $T_pM$ and thereby add them.

Thus the sum on the first line of \eqref{F*Qp} is
\be
\iota_{Y_Q(p)}
\left(\frac{1}{N}\sum_{i=1}^N \iota_{z_i}^{-1}\left[((\exp_p)^{-1})_{*q_i}(w_i)
\right]\right).
\label{F*Qp-1}
\ee
The term on the second line of \eqref{F*Qp} is simply $Y_{Q*p}(v)$. The result now follows
from \eqref{F*Qp}, \eqref{F*Qp-1}, and \eqref{YQ*p}.
\qedns

\begin{cor}\label{Fsubmersive} $F|_{\M_1}$ is a submersion.
\end{cor}

\pf Let $(Q,p)\in \M_1$ and let $u=F(Q,p)$. Let $X\in T_u(TM)$ and let $X_1=\hor_u(X),
X_2=\iota_u^{-1}(\ver_u(X))$, where $\hor_v: T_v(TM)\to \Hor_v(TM)$ and $\ver_v:T_v(TM)\to \Ver_v(TM)$
are the horizontal and vertical projection-maps determined by the splitting of $T_v(TM)$ as $\Hor_v(TM)\plus\Ver_v(TM)$.
Let $v=\pi_{*u}(X_1), z_1=(\exp_p)^{-1}(q_1),
w_1=N(\exp_p)_{*z_1}(\iota_{z_1}(X_2-\na_v Y_Q)),$ and $w_i=0$ for $i>1$.
Then from \eqref{F*qp.0} we have $F_{*(Q,p)}(w_1,0,\dots,0,v)=X_1+X_2=X.$ Hence
$F_{*(Q,p)}$ is surjective. \qedns

Since $Z$ is a codimension-$n$ submanifold
of $TM$, Corollary \ref{Fsubmersive} immediately implies that $(F|_{\M_1})^{-1}(Z)$ is a codimension-$n$ submanifold
of the $(N+1)n$-dimensional manifold $\M_1$. Hence, remembering \eqref{finvz-2}:

\begin{cor}\label{cor:barymfd} $\B_N^0$ is an $(Nn)$-dimensional submanifold of $M^N\times M.$
\qedns
\end{cor}

It is natural to ask how nicely the projection-maps $\pi_N$ and $\plast$ restrict to
$\B_N^0$.  We will see below that $\plast$ restricts to a submersion.
Since $\B_N^0$ and $M^N$ have the same dimension, one may wonder whether
$\pi_N$ restricts to a local diffeomorphism.  We are able to provide only a partial
answer to this question, exhibiting a large open subset of $\B_N^0$ on which $\pi_N$ restricts
this nicely.

Letting $\hess(f_Q)|_p$ denote the covariant Hessian of $f_Q$ at $p$, define
\ben
\M_2 = \{(Q,p)\in \B_N^0  : \hess(f_Q)|_p\ \mbox{is nondegenerate}\},
\een
an open subset of $\B_N^0$. Note that since $Y_Q=-\grad f_Q$, at any $p\in M_Q$
the endomorphism $-(\na Y_Q)|_p$ is the image of $\hess(f_Q)|_p$
under the ``index-lowering'' isomorphism $T_p^*M\tensor T_p^*M\to
T_pM\tensor T_p^*M \caniso \End(T_pM)$ induced by the metric $g$.
Thus the nondegeneracy condition in the definition of $\M_2$ is equivalent to
$(\na Y_Q)|_p : T_pM \to T_p M$ being an isomorphism.

\begin{prop} \label{bary_projected}
The map $\plast|_{_{\B_N^0}}: \B_N^0\to M$ is a submersion, and $\pi_N|_{_{\M_2}}: \M_2\to M^N$ is a local
diffeomorphism.
\end{prop}

\pf Let $(Q,p)=(q_1,\dots,q_N,p)\in \B_N^0$; note that $F(Q,p)=0_p$. Since $T_{0_p}Z$ is exactly the horizontal
subspace of $T_{0_p}(TM)$,
\bestar
\lefteqn{T_{(Q,p)} (\B_N^0)
=(F_{*(Q,p)})^{-1}(T_{0_p}Z)} \\
&=& \{ W\in T_{(Q,p)}M : \ver_{0_p}(F_{*(Q,p)}(W))=0\}\\
&=&
\{ (w_1,\dots, w_N,v)\in T_{q_1}M \times \dots\times T_{q_N}M
\times T_p M : \\
&& \mbox{\hspace{1in}}\na_v Y_Q +
\frac{1}{N}\sum_{i=1}^N \iota_{z_i}^{-1}\left[((\exp_p)^{-1})_{*q_i}(w_i)\right]=0\},
\eestar
\noi where $z_i=(\exp_p)^{-1}(q_i)$; the final equality follows from \eqref{F*qp.0}.
Analogously to the proof of Corollary \ref{Fsubmersive}, let $v\in T_pM$ and set
$w_1=N(\exp_p)_{*z_1}(\iota_{z_1}(-\na_v Y_Q)),$ $w_i=0$ for $i>1$. Then
\be\label{TFinvZ}
\na_v Y_Q +
\frac{1}{N}\sum_{i=1}^N \iota_{z_i}^{-1}\left[((\exp_p)^{-1})_{*q_i}(w_i)\right]=0,
\ee
so $W:=(w_1,\dots,w_N,v)\in T_{(Q,p)}(\B_N^0)$, and $\plast_{*(Q,p)}(W) = v$.  Thus
$\plast_{*(Q,p)}$ is surjective, and $\plast|_{_{\B_N^0}}$ is a submersion.
Now assume that $(Q,p)\in \M_2$ and let $w_i\in T_{q_i}M, 1\leq i\leq N$.
Then the map $(\na Y_Q)|_p : T_pM \to T_p M$ is invertible, and we may define
$$v= -((\na Y_Q)|_p)^{-1}\left(\frac{1}{N}\sum_{i=1}^N \iota_{z_i}^{-1}\left[((\exp_p)^{-1})_{*q_i}(w_i)\right] \right)\in T_p M.$$  Then $W:=(w_1,\dots,w_N,v)\in T_{(Q,p)}(\B_N^0)=T_{(Q,p)}\M_2$, and $\pi_{N*(Q,p)}(W)
= (w_1,\dots,w_N)$.  Thus $(\pi_N|_{\M_2})_{*(Q,p)}$ is surjective.  Since $\dim(\M_2)=\dim(\B_N^0)
=Nn=\dim(M^N)$, \linebreak
$(\pi_N|_{\M_2})_{*(Q,p)}$ is an isomorphism, and $\pi_N|_{_{\M_2}}$ is a local diffeomorphism.
\qedns

\begin{remark}\label{rmk:pi_NM_2}\rm
The result concerning $\M_2$ in Proposition \ref{bary_projected}
is not wholly satisfying, since the set $\M_2$ itself is not easy
to get one's hands on explicitly.  However, at least the set $\pi_N(\M_2)$---the set
of configurations $Q$ for which $\hess(f_Q)|_p$ is nondegenerate for every short barycenter
$p$ of $Q$---is generic in $M^N$.   In  Appendix \ref{sec:appendixB}, we prove this using a
``Parametric Transversality Theorem'' (Theorem \ref{PTT}).
Since $\B_N^0\supset \M_2$, it
follows that $\pi_N(\B_N^0)$ is generic in $M^N$ as well, as asserted earlier.
\end{remark}

\begin{cor}\label{cor:codimk} Let $A\subset M$
be a submanifold of codimension $k$. Then the subset
$\pi_N((\plast|_{_{\M_2}})^{-1}(A))$ of
\ben
\pi_N((\plast|_{_{\B_N^0}})^{-1}(A))=\{Q\in M^N\ : \  \B^0(Q)\intersect A\neq\emptyset\}
\een
is a codimension-$k$ submanifold of $M^N$.
\end{cor}

\pf
Since $\M_2$ is an open subset of $\B_N^0$, and $\plast|_{_{\B_N^0}}$ is a submersion,
$\plast|_{_{\M_2}}$ is also a submersion.  Hence $(\plast|_{_{\M_2}})^{-1}(A)$
is a codimension-$k$ submanifold of $\M_2$, and its image under the local
diffeomorphism $\pi_N|_{_{\M_2}}$ is a codimension-$k$ submanifold of $M^N$.
\qedns

By imposing certain hypotheses under which Fr\'{e}chet means are unique, we will provide in Corollary \ref{cor:conFMbary} a reasonably large and {\em explicit} open subset of $\M_2$
(hence an explicit subset of $\B_N$ that {\em is} a submanifold of $M^N\times M$).
First, some notation:

\begin{defn}\label{def:rinj_rcx}\rm
Let $\rinj=\rinj(M,g)\geq 0$ denote the injectivity radius of $(M,g)$, let
$\Delta=\Delta(M,g) \leq \infty$ be the
the supremum of the sectional curvatures of $(M,g)$, and let
$\rcx=\rcx(M,g)=\frac{1}{2}\min\{\rinj,\frac{\pi}{\sqrt{\Delta}}\}$, where if $\Delta\leq 0$ we interpret $\frac{1}{\sqrt{\Delta}}$ as $\infty$.  Define
\ben
M^N_\con:=\{Q\in M^N:
\supp(Q)
\subset B_{\rcx}(p)\ \mbox{for some $p\in M$}\},
\een
the set of (somewhat) ``concentrated'' configurations.
\end{defn}

Note
that $\rcx>0$ if (and only if) $\rinj>0$ and $\Delta<\infty$, and that $M^N_\con=\emptyset$
if $\rcx=0$. (Hence any hypothesis of the form ``$Q\in M^N_\con$'' implicitly assumes $\rcx>0$.)

A fundamental theorem of Afsari \cite{afsari}, applied to the probability measures $\mu_Q$, yields the following:

\begin{thm}[Special case of {\cite[Theorem 2.1]{afsari}}]\label{thm:afsari}
If $Q\in M^N_\con$, then $Q$ has a unique Fr\'{e}chet mean.
Furthermore, for any open ball $B$ of radius less than $\rcx$
containing $\supp(Q)$, this unique element of $\FM(Q)$ lies in $B$ and is the unique
short barycenter of $Q$ in the concentric open ball of radius $\rcx$.
\end{thm}

\begin{remark}\label{stronger1}
\rm
Theorem 2.1 in \cite{afsari} is far more extensive than Theorem \ref{thm:afsari}, and for digestibility's sake, \
its author makes certain simplifications.  By removing these simplifications, the original theorem and the 
special case above can be strengthened.  For example, instead of the global geometric invariant $\rcx$,  which is
a convenient (and often sharp) lower bound on what is termed the ``regular convexity radius'' of $(M,g)$
in \cite{groisser2004}, local invariants can be used; see \cite[Remark 2.5]{afsari}.
However, because of its simplicity and wide applicability, $\rcx$ is very commonly used in hypotheses on
the radius of a ball supporting a probability distribution for which one wants to guarantee
the uniqueness of a Fr\'{e}chet mean.
\end{remark}

Let $Q\in M^N_\con$.
The proof of Theorem 2.1 in \cite{afsari} shows that
$\hess(f_Q)$ is positive-definite, and hence nondegenerate, at the (unique) Fr\'{e}chet mean of $Q$.
As a corollary, we have:

\begin{cor} \label{cor:conFMbary} The set
\be\label{conFMbary.1}
\M_{\rm con}:=\{(Q,p)\in \B_N \ : \  Q\in M^N_\con\}= (\pi_N|_{\B_N})^{-1}(M^N_\con)
\ee
is a submanifold of $M^N\times M$, and the restriction of $\pi_N$ to this submanifold
is a diffeomorphism onto its image.
\end{cor}

\pf
Theorem \ref{thm:afsari} implies that if $Q\in M^N_\con$ and $p$ is the Fr\'{e}chet mean of $Q$, then
$(Q,p)\in \M_2$.  Thus $\M_{\rm con}$ is a subset
of the manifold $\M_2$, and is open in $\M_2$ (hence is a submanifold of $M^N\times M$) since $M^N_\con$ is open in $M^N$.  The map $\pi_N|_{_{\M_{\rm con}}}$ is the restriction of the local
diffeomorphism $\pi_N|_{_{\M_2}}$ to the open subset $\M_{\rm con}\subset \M_2$, and
is one-to-one by Theorem \ref{thm:afsari}.
\qedns

\begin{remark}\label{rem:findiff}
 \rm As noted in the introduction,  we have taken ``smooth'' to mean $C^\infty$, and have taken
``(sub)manifold'' to mean ``smooth (sub)manifold''.  However, all statements in Sections \ref{sec:bary-prebary2} and 
\ref{sec:barymanifold} up through Corollary \ref{cor:codimk} remain true (and the same proofs  work) if the manifold $M$ and Riemannian metric $g$ are 
assumed only to be of class $C^k$, where $k\geq 2$, and other explicit or implicit references to smoothness are modified
as follows: 

\begin{itemize}
\item The functions $f_Q$ are $C^k$.  

\item The tangent bundle $TM$ is a $C^{k-1}$ manifold.  The maps $\exp, \exp_p,\pi, \pi_N,\plast,$ and $F_Q$ are $C^{k-1}$, as are the vector fields $Y_Q$ and the submanifolds $Z\subset TM$, $\B_N=F^{-1}(Z)\subset  M^N\times M$, and $\M_2\subset M^N\times M$.  

\item The maps $\plast|_{_{F^{-1}(Z)}}$ and $\pi_N|_{_{\M_2}}$ are $C^{k-1}$.

\item The submanifold $\pi_N((\plast|_{_{\M_2}})^{-1}(A))$ in Corollary \ref{cor:codimk}
is $C^1$.

\end{itemize}

\noi However, for the genericity statement in Remark \ref{rmk:pi_NM_2}, a much 
higher degree of smoothness of $(M,g)$ is needed:  we need to assume that $M$ and $g$ are $C^k$, where 
$k\geq n(N-1)+2.$ The reason is that,
just as in the finite-differentiability version of Sard's Theorem \cite[Theorem III.1.3]{hirsch} of which 
the finite-differentiability version of the Parametric Transversality Theorem \cite[Theorem III.2.7]{hirsch}
is a corollary,  the required degree of smoothness depends on the dimensions of the manifolds involved, via an
explicit inequality that reduces to the stated restriction on $k$ in our setting.
\end{remark}

\setcounter{equation}{0}
\section{All Fr\'{e}chet means of finite configurations are barycenters}
\label{sec:FMs_are_barys}

Theorem \ref{thm:afsari} shows that the Fr\'{e}chet mean set of a sufficiently
concentrated configuration $Q\in M^N$ lies in $\B(Q)$.
Using results from \cite{le-barden2014}, we will remove the ``sufficiently concentrated'' restriction.

\ss
For any $f:M\to\bfr, p\in M,$ and $v\in T_pM$, define the {\em forward geodesic directional derivative}
$v_+(f)$ by
\be\label{fdirder}
v_+(f)=\lim_{t\downto 0} \frac{f(\exp_p(tv))-f(p)}{t}
\ee
whenever the limit exists.
Note that if $f$ has a relative minimum at $p$, and $v_+(f)$ exists, then
$v_+(f)\geq 0$.

\begin{remark}\label{rem1:direcderiv} \rm
In \cite[Section 11]{burago}, in the more general setting of Alexandrov spaces with curvature bounded below,
the notion corresponding to \eqref{fdirder}  is developed (only) for functions $f$ {\em that are locally Lipschitz at $p$}---i.e.,
that satisfy a Lipschitz condition on some neighborhood of $p$.
{\em With this restriction on $f$}, the limit \eqref{fdirder} is simply called a ``directional
derivative'' in \cite{burago}.  This terminology, and the notation $v(f)$ for our $v_+(f)$,
are used in \cite{le-barden2014} for powers (with exponent $\geq 1$) of the distance function $r_p=d(p,\cdot).$ 
However, the omission from \cite{le-barden2014} of any mention of Lipschitzness may leave a misperception that 
the notation is unambiguous, and consequently that certain facts discussed below are true,
for any continuous function $f:M\to\bfr$. We use the longer ``forward geodesic directional derivative''
because of several sources of potential misunderstanding, all related to the fact that
in the setting of smooth manifolds, there is a conventional---and not metric-related---use of the term ``directional derivative''
that, unlike in vector spaces, is usually applied only to {\em differentiable}
functions.  If $f:M\to\bfr$ is differentiable at $p$ (below we leave the ``at $p$'' implicit), and $v\in T_pM$, then what is commonly
called the ``directional derivative $v(f)$'' is the two-sided limit $\frac{d}{dt}f(\g(t))|_{t=0},$ where $\g:I\to \bfr$ is {\em any}
smooth curve based at $p$ with $\g'(0)=v$. Crucially, differentiability of $f$ ensures
that (i) $\frac{d}{dt}f(\g(t))|_{t=0}$  is independent of the
choice of curve $\g$ representing $v$, so that the notation ``$v(f)$'' makes sense, and
(ii) $v_+(f)=v(f)$, so that the word ``forward'' would be superfluous.
But if $f$ is not differentiable at $p$, then even if $(-v)_+(f)=-(v_+(f))$ (so that $\frac{d}{dt}f(\exp_p(tv))|_{t=0}$ exists), it is not always true that for {\em every} curve $\g$ representing $v$, the derivative $\frac{d}{dt}f(\g(t))|_{t=0}$  exists or has 
the same value. Even the existence of the map $(d_+f)_p: T_pM\to\bfr$, $v\mapsto v_+(f)$ (i.e. existence of
$v_+(f)$ for all $v\in T_pM$) plus linearity in $v$, are not sufficient to ensure the curve-independence property.  
However, for functions that are locally
Lipschitz at $p$, existence  {\em and continuity} of the map $(d_+f)_p$ together ensure the curve-independence property.  
(This can be deduced from the arguments in \cite{burago}.)  Therefore for
positive powers of the distance function $r_p$, the terminology ``directional derivative'' and notation 
$v(f)$ in \cite{le-barden2014} are consistent with the usual terminology and notation for functions on
manifolds, modulo the one-sidedness of the limit \eqref{fdirder}.
\end{remark}

For $q\in M$, the function $r_q^2:M\to\bfr$ is never differentiable at ordinary cut-points of $q$ \cite{bishop}.
However, Le and Barden show in \cite{le-barden2014} that at {\em every} cut-point of $q$, all forward geodesic directional derivatives of $r_q^2$ exist. An explicit formula, which is valid also at non-cut points of $q$, is derived there:

\begin{lemma}[{\cite[Lemma 2]{le-barden2014}}]
\label{bllemma2}
Let $p,q\in M$,
$v\in T_pM$.  Then
\be
v_+(r_q^2) = -2\sup_{v'\in \V_{q,p}} g(v,v'),
\ee
where
$\V_{q,p}=\{v'\in T_pM \ : \  \|v'\|=d(q,p)\ \mbox{\rm and}\ \exp_p(v')=q\}.$
\end{lemma}

As noted in \cite{le-barden2014}, if $p\in\clord(q)$, then since $r_q^2$ is not differentiable at $p$, we do not necessarily have $(-v)_+(r_q^2)$ equal to $-(v_+(r_q^2))$.  Also note that if $p\notin \clord(q)$, then $\V_{q,p}$ consists of the single element $(\exp_p)^{-1}(q)$.

The proofs of Theorem 1 and Corollaries 2 and  3 in \cite{le-barden2014}
show the following:  Let $\mu$ be probability measure on $M$ that
is the sum of an absolutely continuous measure and a discrete measure, let $\a>1$, and
assume that the $\a$-energy function
$f_{\mu,\a}: x\mapsto \frac{1}{\a}\int_M d(x,y)^\a \, d\mu(y)$ achieves a relative minimum at $p$. Then (i) if
$\clsing(p)=\emptyset$, then $\mu(\CL(p))=0$;  (ii) $\exp_p^{-1}(q)$ is defined
$\mu$-a.e. in $q$; and (iii)
$\int_M \exp_p^{-1}(q) d\mu(q)=0.$  Proving fact (i) was the main goal of \cite{le-barden2014};
for us, more directly important
is fact (iii) (but note that the integral in (iii) is well-defined only because of (ii)).
Specializing these facts to the case
$\mu_Q=\frac{1}{N}\sum_i \d_{q_i}$ and $\a=2$ yields Proposition \ref{prop:fmisbary}
below.  We will prove this proposition more directly, however.  Our proof borrows heavily from the proof of \cite[Theorem 1]{le-barden2014}, but our argument is rather simpler.

\begin{prop}\label{prop:fmisbary}  Let $Q=(q_1,\dots,q_N)\in M^N,$ let $p\in M$ and assume that $f_Q$ achieves a relative
minimum at $p$. Then:

\begin{itemize}
\item[(i)] For $1\leq i\leq N$ we have $q_i\in \tD(p)$, and hence $\exp_p^{-1}(q_i)$
is defined.

\item[(ii)] The point $p$ is a barycenter of $Q$ (in fact, an almost-short barycenter).
\end{itemize}

\end{prop}

\pf For (i), it suffices to show that $\supp(Q)\intersect \clord(p)=\emptyset$. Let $I_{\rm reg}, I_\ord,$ and  $I_\sing$ denote the subsets of $\{1,\dots, N\}$ for which
$q_i$ lies, respectively, in $\D(p), \clord(p)$, and $\clsing(p)$.  Then for every $v\in T_pM$,
by Lemma \ref{bllemma2} we have
\be\label{v+fq}
v_+(f_Q)= - \frac{1}{N}\left\{\sum_{i\in I_{\rm reg}\union I_\sing}
g(v,\exp_p^{-1}(q_i))
+\sum_{i\in I_\ord}
\sup_{v'\in \V_{q_i,p}} g(v,v')\right\}.
\ee
Since $f_Q$ has a relative minimum at $p$, we have $v_+(f_Q)\geq 0$ for every $v\in T_pM$.  Hence,
for each $v\in T_pM$,
\bearray\nonumber
0&\leq &
N[v_+(f_Q)+ (-v)_+(f_Q)]\\
\nonumber
&=&
-\sum_{i\in I_\ord}\left( \sup_{v'\in \V_{q_i,p}} g(v,v')
+\sup_{v'\in \V_{q_i,p}} g(-v,v')\right)\\
&=&
\sum_{i\in I_\ord}\left( \inf_{v'\in \V_{q_i,p}} g(v,v')
-\sup_{v'\in \V_{q_i,p}} g(v,v')\right)
\label{inf-sup}
\\
&\leq & 0,
\nonumber
\eearray
implying that each summand in \eqref{inf-sup} is zero.
Hence for each $i\in I_\ord$ and $v\in T_pM$, the set $\{g(v,v')\ : \  v'\in \V_{q_i,p}\}$
consists of a single element, implying that the difference of any two elements of
$V_{q_i,p}$ is perpendicular to every element of $T_pM$, hence is zero.
But by definition of $\clord(p)$, if $q\in \clord(p)$ then $V_{q,p}$ contains
at least two vectors.  Hence $I_\ord=\emptyset$, establishing (i).

Substituting $I_\ord=\emptyset$ into \eqref{v+fq}, we now have
$v_+(f_Q)= - \frac{1}{N}\sum_{i=1}^N
g(v,\exp_p^{-1}(q_i))\geq 0$
for all $v\in T_pM$.  Replacing $v$ with $-v$, we obtain the opposite inequality,
and hence the equality $\sum_{i=1}^N
g(v,\exp_p^{-1}(q_i))=0$.  Thus $p\in \tB^0(Q)$.\qedns

\begin{remark}\label{rem2:direcderiv}  \rm
The proof above also shows that $v_+(f_Q)=0$ for all $v\in T_pM$.  For a general
function $f$, having $v_+(f)=0$ for all $v\in T_pM$ does {\em not} imply that $f$ is differentiable at
$p$;
see Remark \ref{rem1:direcderiv}.  However, it is not hard to show that $f_Q$ is locally Lipschitz
(or that the same is true more generally for the function $f_{\mu,\a}: p\mapsto 
\int_M r_p^\a\, d\mu$,
where $\a\geq 1$ and $\mu$ is any probability measure on $M$ for which $f_{\mu,\a}$ is finite).
Therefore for $f_Q$, the fact that $v_+(f_Q)=0$ for all $v\in T_pM$ {\em does} imply that
$f_Q$ is differentiable at $p$;
again see Remark \ref{rem1:direcderiv}.
\end{remark}

An immediate corollary of Proposition \ref{prop:fmisbary} is the following:

\begin{cor}\label{cor:fmisbary}
For every $Q\in M^N$, we have $\FM(Q)\subset \tB^0(Q)\subset \B(Q)$.\qedns
\end{cor}

\noi {\bf Proof of Theorem \ref{thm:FMmeaszero}.} The theorem follows immediately from Corollary \ref{cor:fmisbary} and  Theorem \ref{thm:measzero.7}.
\qedns

\setcounter{equation}{0}
\section{Fr\'{e}chet sample means of i.i.d. random variables}
\label{sec:iid1}

$~$

\ss
{\bf Convention.}
 For the remainder of this paper, given a probability measure $\mu$
on a topological space $Y$, and a subset $A\subset Y$, we will say that {\em $\mu$ is supported
in $A$} if $\mu(A)=1$.  This is less restrictive than saying that $A$ contains the set $\supp(\mu)$
(the support of $\mu$),  which by definition is a closed set.

\begin{cor}\label{cor:measzero.3}
Let
$X_1, \dots, X_N$ be absolutely continuous,
i.i.d.
$M$-valued variables with underlying measure $\mu_X$,
and let $A$ be a volume-zero subset of $M$.  Then:
\begin{enumerate}
\item[(a)]
${\rm Pr}\left\{\FM(X_1,\dots,X_N)\subset M\backslash A\right\}=1.$

\item[(b)]
If  $\mu_X$ is supported
in  an open ball $B$ of radius $r\leq\rcx$, then almost surely the Fr\'{e}chet mean set of $(X_1,\dots,X_N)$
is a single point and lies in $B\backslash A$.
\end{enumerate}
\end{cor}

\pf Let $\mathcal{X}$ denote the $M^N$-valued random variable $(X_1,\dots,X_N)$
and let $\mu_X^N$ be the product
probability measure on $M^N$ induced by $\mu_X$. Since $\mu_X$ is absolutely continuous,
so is $\mu_X^N$.
Hence
\ben
{\rm Pr}\left\{\FM(\mathcal{X})\intersect A\neq \emptyset\right\}
= \mu_X^N(\{Q\in M^N: \FM(Q)\intersect A\neq \emptyset\})
=0
\een
by  Theorem \ref{thm:FMmeaszero}.  Statement (a) follows.

Now suppose that $\mu_X$ is supported in an open ball $B$ of radius $r\leq\rcx$, and let $B^N$ denote the $N$-fold 
Cartesian product of $B$ with itself. Then ${\rm Pr}\{\mathcal{X}\in B^N\}=1$.  If the support of a configuration 
$Q\in M^N$ lies in the open ball $B$, then $\supp(Q)$ also lies
in a concentric open ball of  radius strictly less than $r$, hence less than $\rcx$. Therefore Theorem \ref{thm:afsari}
ensures that if a value of $\mathcal{X}$ lies in $B^N$, then the Fr\'{e}chet mean
of that value of $\mathcal{X}$ is unique and lies in $B$.  Statement (b) follows.
\qedns

\begin{remark}\label{weakstrat}
\rm
Corollary \ref{cor:measzero.3}(a) is precisely Theorem \ref{thm1.1}. This theorem has an application to manifolds 
that are stratified in even a very weak sense.  Let us call $M$ {\em weakly stratified} if we are given a collection 
$\mathcal T$ of pairwise disjoint submanifolds of $M$, the {\em strata}, whose union is $M$, and for
which exactly one stratum, $M_{\rm top}$, is an open set. (More precisely, the ``weakly
stratified manifold'' is the pair $(M,\mathcal{T})$, but in many situations,  it is often not
important to have a symbol for the set $\mathcal T$.) Then each of the other
strata has positive codimension, so the union of all the other strata has volume zero in $M$.  Hence Theorem \ref{thm1.1}
implies that if $X_1, \dots, X_N$ are absolutely continuous,  i.i.d., $M$-valued random variables, then
${\rm Pr}\left\{\FM(X_1,\dots,X_N)\subset M_{\rm top}\right\} =1.$  

But note that Theorem \ref{thm1.1} is much more general than 
the above application; the theorem does not require {\em any} sort of stratification.
\end{remark}

An example of a well known manifold to which Theorem \ref{thm1.1} applies is the projective space ${\bf C}P^{k-2}$.  This manifold arises in the statistics of shape as Kendall's {\em shape space}  $\Sigma^k_2$ of $k>2$ labeled points in $\bfr^2$
(``planar shapes'') \cite{kendall_d}.
With the usual shape-space stratification,  $\Sigma^k_2$ has only one stratum (the entire space),
so the results of \cite{huckMMS2012} reduce to tautologies.   However,  $\Sigma^k_2$ has other
stratifications.  For example,
we can stratify according to the greatest number of collinear points in the configuration.  In this
case the top stratum consists of the shapes whose pre-shapes lie in {\em general position}
(no three points are collinear), and the remaining strata have positive codimension.   Theorem
\ref{thm1.1} shows that for i.i.d. random planar shapes $X_1,\dots, X_N$ with absolutely continuous
distribution, the sample Fr\'{e}chet mean shape almost surely has no more than two collinear points.
More generally, a similar result holds for many other ways in which a planar shape can be non-generic.
For example, for $k\geq 4$, if we call a point in $\Sigma^k_2$ a {\em $k$-trapezoid} if it has $k$ distinct points and at least two parallel non-adjacent sides, then
the set $\Sigma^{k, {\rm trap}}_2$ of $k$-trapezoids is a finite union of codimension-1 submanifolds of $\Sigma^k_2$.  Thus if $X_1,\dots, X_N$ are as above, then by Theorem \ref{thm1.1},
${\rm Pr}\left\{\FM(\{X_1,\dots,X_N\})\intersect \Sigma^{k,{\rm trap}}_2\neq\emptyset\right\}=0$.

Theorem \ref{thm1.1} applies only to  {\em manifolds}, but allows {\em arbitrary} (volume-)measure-zero subsets $A$.  There are certain, more general, {\em non-manifold} spaces for which similar results have been obtained, but
only for {\em very specific} measure-zero subsets $A$:  unions of lower strata of an orbit-type stratification
of a quotient space $M/G$, where $M$ is a complete, connected Riemannian manifold and $G$ is a Lie group acting properly and isometrically on $M$.  These spaces are (generalized) shape spaces, as broadly defined
in \cite{huckMMS2012}.  The orbit-type stratification
of $M$ induces a stratification of $M/G$, which we will call the ``usual shape-space stratification'', and
$M/G$ inherits a quotient distance-function.  For these spaces the behavior of
sample means is similar to what Theorem \ref{thm1.1} shows what it is for a manifold:  in \cite[Corollary 2]{huckMMS2012},
Huckemann establishes that  Fr\'{e}chet sample means of an absolutely
continuously distributed random variable lies in the top stratum---the complement of
the union $A$ of all lower strata (a particular measure-zero set)---with probability 1.

There are other situations  in which Theorem  \ref{thm1.1} does not apply.
In recent years there has been considerable interest in the behavior of various means on
certain stratified spaces that usually are not manifolds, for example ``open
books'' and ``spiders'' 
%\citep{hotz2013sticky}. 
\cite{hotz2013sticky}. For open-books and spiders, Fr\'{e}chet sample means are 
often ``sticky'', tending to lie in lower strata 
%\citep{hotz2013sticky}.
\cite{hotz2013sticky}. See also 
%\citep{MR3371437}.
\cite{MR3371437}.
The ``sticky means'' behavior on these non-manifold spaces is in sharp contrast to what Theorem \ref{thm1.1}
implies for (even weakly) stratified manifolds.

\setcounter{equation}{0}
\section{Riemannian covers and equivariant means}
\label{sec:eqvt_means_general}

In some applications, the easiest way to compute or analyze Fr\'{e}chet sample
means of configurations (or random variables) is to relate them to what we will
call {\em equivariant Fr\'{e}chet means} on a Riemannian covering space.  We
will see an example of this in Section \ref{sec:PSRmeans}.

Let $G$ be a finite group  of order at least 2,
 with identity element $e$, acting isometrically on
a Riemannian manifold $(\tim,\tg)$; without loss of generality we
assume that the action is a left-action, and denote it by $(\tp,h)\mapsto h\dotprod \tp$. Assume that
$\b:=\inf\{\tilde{d}(\tp,h\dotprod \tp): \tp\in \tim, \ e\neq h\in G\}>0;$
in particular this implies that the $G$-action is free.
For simplicity's sake, we assume throughout this section that $(\tim,\tg)$
is complete, but this assumption is not actually needed until Lemma \ref{lem:rinjbeta}; in particular,
Propositions \ref{quodist} and \ref{prop:uniqueuptoG} do not require completeness.

\begin{remark}\label{princ_rc}\rm  Because the $G$-action is free, $\tim$ is a (left) principal bundle
over $M$ (more strongly, a ``principal Riemannian cover'').  The principality assumption is primarily a convenience, 
allowing simplification of certain statements whose full generality is not needed  for our  application to $\Sym^+(m)$ mentioned in the Introduction.  Many of the 
ideas and results in this section generalize to non-principal Riemannian covers. 
\end{remark}

Since the finite $G$-action is free and isometric, $(\tim,\tg)$ is a {\em Riemannian
covering space}: the quotient $M=\tim/G$ is (canonically) a manifold, the quotient map $\quo:\tim\to M$ is
a covering map, and there is a (unique) induced Riemannian metric $g$ on $M$ such that $\quo$ is a
local isometry.   The assumed completeness of $(\tim, \tg)$
implies that $(M,g)$ is complete as well. We will write $\td$ and $d$ for the Riemannian distance functions on
$(\tim,\tg)$ and $(M,g)$, respectively; in this section, balls are always taken with respect to these distance functions.  
We also use $\td(\cdot,\cdot)$ and $d(\cdot,\cdot)$ to denote the corresponding distances between {\em subsets} 
of $\tim$ and $M$, respectively.

\begin{prop}\label{quodist} For all $p,q\in M$ and $\tp\in\quo^{-1}(p)$,
we have
\be\label{dbar}
d(p,q)=\td(\quo^{-1}(p), \quo^{-1}(q))
=\td(\tp,\quo^{-1}(q)).
\ee
\end{prop}

\pf This is straightforward and is left to the reader.\qedns

The equivariant means defined below are a special case of ``$\rho$-means'' as extended by Huckemann \cite{huck3DProcrustes2011,Huckemann2011consistency}. 

\begin{defn}\label{lem:uniqueuptoG} \rm Define  an ``equivariant distance function'' $\devt:M\times \tim
\to\bfr$  by 
\be\label{defdevt}
\devt(q,\tp)=\td(\quo^{-1}(q),\tp)=\min\{\td(\tq,\tp): \tq\in \quo^{-1}(q)\}.
\ee
For $Q=(q_1,\dots,q_N)\in M^N$ define the ``equivariant Fr\'{e}chet function'' $\tf_{Q}:\tim\to \bfr$ by
\be\label{deftfq}
\tf_{Q}(\tp)=\frac{1}{N}\sum_{i=1}^N \devt(q_i,\tp)^2.
\ee
 Call $\tp\in \tim$ an {\em equivariant Fr\'{e}chet mean} of $Q$
if $\tp$ is a minimizer of $\tf_{Q}$, and write $\EFM(Q)$
for the set of equivariant Fr\'{e}chet means of $Q$.
\end{defn}

Note that by Proposition \ref{quodist}, for all $q\in M$ and $\tp\in \tim$ we have
\be\label{dpar=dbar}
\devt(q,\tp)=d(q,\quo(\tp)),
\ee
and hence
\be\label{fbar=f}
\tf_{Q}(\tp)
=f_{Q}(\quo(\tp))=\tf_{Q}(h\,\dotprod \tp)\ \ \ \mbox{for all $h\in G$}.
\ee
Hence if $\tp$ is an equivariant Fr\'{e}chet mean of $Q$, then so is $h\dotprod \tp$ for every
$h\in G$.
Analogously to \cite{huck3DProcrustes2011} and \cite{rjgs},
we
say that {\em the equivariant Fr\'{e}chet
mean of $Q$ is unique up to the action of $G$} if any two equivariant Fr\'{e}chet
means of $Q$ are related by the action of an element of $G$.

As one might expect, taking equivariant Fr\'{e}chet means in $\tim$ is ``equivalent''
to taking
Fr\'{e}chet means in $M$:

\begin{prop}\label{prop:uniqueuptoG}
Let $Q\in M^N$.
\begin{itemize}
\item[(a)]  Let $\tp\in \tim$ and let $p=\quo(\tp)$. Then the following are equivalent:

\begin{itemize}
\item[(i)] $p\in \FM(Q).$

\item[(ii)] $\tp\in \EFM(Q)$.

\item[(iii)] $\quo^{-1}(p) \subset \EFM(Q)$.

\end{itemize}

\ss\mbox{\hspace{-.2in}} Thus $\EFM(Q)=\quo^{-1}(\FM(Q))$
\  and \ $\quo(\EFM(Q))=\FM(Q)$.

\ss
\item[(b)] $Q$ has a unique Fr\'{e}chet mean if and only if the
equivariant Fr\'{e}chet
mean of $Q$ is unique up to the action of $G$.
\end{itemize}
\end{prop}

\pf This follows from equation \eqref{fbar=f}.\qedns

\begin{cor}\label{EFM_existence}
 Equivariant Fr\'{e}chet means always exist. More precisely, for any \lnbrk 
$Q\in M^N$, the set $\EFM(Q)$ is nonempty.
\end{cor}

\pf  Let $Q\in M^N$. As noted in the introduction, since $(M,g)$ is complete, $\FM(Q)$
is nonempty.  Since $\quo:\tim\to M$ is surjective, $\quo^{-1}(\FM(Q))$ ($=\EFM(Q)$ by 
Proposition \ref{prop:uniqueuptoG})
is nonempty as well. \qedns

To obtain a useful equivariant-means  analog of Corollary \ref{cor:measzero.3},
we need to know how the numbers $\rcx(\tim,\tg), \rcx(M,g),$ and $\b$ are related.  We
start with the following lemma. 

\begin{lemma}\label{lem:rinjbeta}
$\rinj(\tim,\tg)\ \geq\ \rinj(M,g)\ = \ \min\left\{\rinj(\tim,\tg), \frac{\b}{2}\right\}$.
\end{lemma}
This lemma is almost certainly known, but the authors have found no reference for it. For a proof, see
Appendix  \ref{sec:appendixA}.

\begin{cor}\label{cor:rcxbeta} 
$\rcx(\tim,\tg)\ \geq\ \rcx(M,g)\ =\ \min\{\rcx(\tim,\tg),\frac{\b}{4}\}$.
 
\end{cor}

\pf Since $(\tim,\tg)$ is a Riemannian covering space of $(M,g)$, the
sectional-curvature function of $(M,g)$ has the same range as the sectional-curvature function
of $(\tim,\tg)$. Hence $\Delta(M,g)=\Delta(\tim,\tg)$.  The result now follows
from Lemma \ref{lem:rinjbeta}.
\qedns

Since $\quo: \tim\to M$ is a covering map and a local diffeomorphism,
every sufficiently small open set $A\subset M$ is {\em smoothly evenly covered}:
$\quo^{-1}(A)$ is a disjoint union
of open sets $\tilde{A}_j$ (unique up to labeling), each of which is carried diffeomorphically to $A$ by $\quo$.
For any such $A, \{\tilde{A}_j\}$, and arbitrary subset $A'\subset A$, the set $\quo^{-1}(A')$
is the disjoint union of the sets $\tilde{A}'_j:=\tilde{A}_j\cap\quo^{-1}(A')$, each of which is
carried homeomorphically to $A$ by $\quo$. We will
say that $A'$ is evenly covered by the sets $\tilde{A}'_j$ (whether or not $A'$ is open),
and that the collection $\{\tilde{A}'_j\}$ is an {\em even covering} of $A'$ (a {\em smooth even
covering} if $A'$ is open).

\begin{cor}\label{evenlycov}
Let $p\in M$, let $r<\rinj(M,g)$, and let $B=B_r(p)$. Then $B$ is open (as a subset
of the manifold $M$) and is smoothly evenly covered by a collection of $\td$-open
balls $\tilde{B}_j=B_r(\tp_j)\subset (\tim,\td)$, $1\leq j\leq |G|$.  For each $j$, the map $\quo|_{\tilde{B}_j}: \tilde{B}_j\to B$ 
is an isometry.
\end{cor}

\pf Since $d$ is the geodesic-distance function of a Riemannian metric on $M$, the metric-space topology of
$(M,d)$ coincides with the manifold topology. Hence the metrically open
ball $B$ is also open in the manifold topology. (This would {\em not} be true
for an {\em arbitrary} distance-function on $M$.)

Since $r<\rinj(M,g)\leq \rinj(\tim,\tg)$, the exponential maps $\exp_{p}$, $\exp_\tp$,
restricted to the open balls of radius $r$ centered at $0$ in $T_{p}M$, $T_\tp\tim$
respectively, are diffeomorphisms onto $B_r(p)$, $B_r(\tp)$, respectively.
Equation \eqref{intertwine}  in the proof of Lemma \ref{lem:rinjbeta}
 then implies that $\quo|_{B_r(\tp)}$ is a diffeomorphism
onto $B_r(p)$.  Since $\quo$ is a local isometry, it follows that $\quo|_{B_r(\tp)}$
is an isometry.

It remains only to show that for distinct $\tp_1,\tp_2\in \quo^{-1}(p)$, the sets
$\tilde{B}_j:=B_r(\tp_j)$ ($i=1,2$) are disjoint.  Suppose there exists $q\in
\tilde{B}_1\cap\tilde{B}_2$.  Then there are minimal geodesics $\g_j:[0,1]\to \tim$ from $\tp_j$
to $\tq$ ($i=1,2$), each of length less than $r$, and $\g_1'(1)\neq \g_2'(1)$
(else both $\tp_1$ and $\tp_2$ would be images of the same vector in $T_\tq\tim$
under $\exp_\tq$).
Hence $\bgam_j:=\quo\circ\g_j$ ($j=1,2$)
is a geodesic from $p$ to $q\in B_r(p)$ of length less than $r$,
and $\bgam_1'(1)\neq \bgam_2'(1)$ (since $\quo_{*\tq}$ is an isomorphism);
thus $\bgam_1\neq\bgam_2$.  But then the restriction of $\exp_{p}$ to
the ball $B_r(0_{p})\subset T_{p}\tim$ is not one-to-one, a contradiction.
Hence $\tilde{B}_1\cap\tilde{B}_2=\emptyset.$
\qedns

\begin{remark}\label{relabel}\rm
In the setting of Corollary \ref{evenlycov}, the $G$-action permutes the sets $\tilde{B}_j$.
While there is, in general, no {\em canonical} labeling of the $\tilde{B}_j$
by elements of $G$, if we arbitrarily label any one of these balls as $\tilde{B}_e$, then we can
(re-)label all of the others by setting $\tilde{B}_h=h\,\dotprod \tilde{B}_e$ for all $h\in G$.  We then
have $h_1\,\dotprod \tilde{B}_{h_2}=\tilde{B}_{h_1h_2}$ for all $h_1,h_2\in G,$ so
we call such a labeling {\em equivariant}.
This facilitates the statement of the next proposition, a corollary of earlier results.
\end{remark}

\begin{prop}\label{prop:equivar_means} Let ${r}\in (0,\rcx(M,g)]$ and
let $Q=(q_1,\dots,q_N)\in M^N$.  Assume that
$\supp(Q)$ lies in a ball $B$ of radius less than $r$.
Let
$\{\tilde{B}_h\}_{h\in G}$
be the
even covering of $B$ given by Corollary \ref{evenlycov},
with the covering balls labeled equivariantly,
and for each $h\in G$ and $p\in B$,
let $\tilde{p}^{(h)}$
be the unique element of $\quo^{-1}(p)\cap \tilde{B}_h$.
Let $\tQ^{(h)}=(\tilde{q}_1^{(h)}, \dots, \tilde{q}_N^{(h)})$. Then

\begin{itemize}
\item[(a)] For each $h\in G$, $\FM(\tQ^{(h)})$ consists of a single point
and lies in $\tilde{B}_h$.

\item[(b)] For all $h_1,h_2\in G$ we have $h_1\dotprod \FM(\tQ^{(h_2)})
=\FM(\tQ^{(h_1h_2)}).$

\item[(c)] $\EFM(Q)=\{\FM(\tQ^{(h)})\}_{h\in G}
\subset \quo^{-1}(B').$  In particular, the equivariant Fr\'{e}chet
mean of $Q$ is unique up to the action of $G$.

\item[(d)] $\FM(Q)=\quo(\EFM(Q))=\quo(\FM(\tQ^{(h)}))$ for all $h\in G$,
and consists of a single point lying in $B$. Since
$\EFM(Q)$ is a $G$-invariant set, the
first equality is equivalent to
$\EFM(Q)=\quo^{-1}(\FM(Q))$.
\end{itemize}

\end{prop}

\pf (a) Let $h\in G$.
By Corollary \ref{evenlycov}, $\tilde{B}_h$ is a ball of radius $r<\rcx(M,g)\leq \rcx(\tim,\tg)$.
Hence
Theorem \ref{thm:afsari} implies that the Fr\'{e}chet mean of $\tQ^{(h)}$ is unique
and
lies in $\tilde{B}_h$.

\ss (b) Since $h_1\dotprod B'_{h_2}= B'_{h_1h_2}$ we have
$h_1\dotprod \tQ^{(h_2)}= \tQ^{(h_1h_2)}$. Since the metric on $M$ is $G$-invariant,
it follows that $h_1\dotprod \FM(\tQ^{(h_2)})= \FM(\tQ^{(h_1h_2)})$.

\ss (c) Follows from Proposition \ref{prop:uniqueuptoG}(a).

\ss (d) Follows from (c), (a), and
Proposition \ref{prop:uniqueuptoG}(b).
\qedns

\setcounter{equation}{0}
\section{Fr\'{e}chet sample means in the equivariant setting}
\label{sec:iid_eqvt_means}

$~$

In this section, $(\tim,\tg)$, $G$, $(M,g)$, and $\quo$ are as in Section \ref{sec:eqvt_means_general}.

Given any set $S$, any function $f:S\to \tim$, and any $h\in G$, let $h\dotprod f$ denote the function
$x\mapsto h\dotprod f(x)$.  This defines a (left) action of $G$ on the set of functions from $S$ to $\tim$. 
We will use the term {\em minimal invariant family of functions (from $S$ to $\tim$)} for an orbit of $G$ under this action, and 
{\em equivariant family of functions  (from $S$ to $\tim$)} for an indexed collection $\{f^{(h)}:S\to \tim\}_{h\in G}$ satisfying
$f^{(h_1h_2)}=h_1\dotprod f^{(h_2)}$ for all $h_1,h_2\in G$.  (In the latter case, we are opting
for terminology that is simpler than the more precise ``{\em equivariantly indexed minimal family of functions  (from $S$ to $\tim$)}''.) For these families of functions, there will be only two types of domains of concern to us: the domain $\Omega$ of a random variable, and subsets of $M$.  In each case the domain $S$ will be clear from context, so we will usually omit the ``from $S$ to $\tim$.'' 

Clearly the underlying set of functions in an equivariant family of functions is a minimal invariant family.  Conversely, given a minimal invariant family $\F$ , any chosen element $f_1\in \F$ determines an equivariant indexing of $\F$ by setting $f{(h)}=h\dotprod f_1$ for each $h\in G$.

Corollary \ref{evenlycov} and Remark \ref{relabel} illustrate one way that such families of functions arise naturally. Consider a $d$-open ball $B\subset M$ that is smoothly evenly covered by a collection of $\td$-open balls $\tilde{B}_j$, $1\leq j\leq |G|$.
For each such $j$ let $s_j=(\quo|_{\tilde{B}_j})^{-1}:B\to \tilde{B}_j$, but viewed as a map from $B$ to $\tim$ (a {\em local section} of the bundle $\tim\stackrel{\quo}{\longrightarrow} M$) . Then 
$\{s_j: 1\leq j\leq |G|\}$ is a minimal invariant family of smooth maps.  A re-indexing of the collection $\{\tilde{B}_j\}$ by the elements of $G$, as in Remark \ref{relabel}, yields a corresponding re-indexing of this collection of local sections.  The re-indexed family $\{s_h\}_{h\in G}$ is an equivariant family of smooth maps.

\begin{defn} \rm 
We will call a map $s:M\to\tim$ a {\em measurable section} (of the $G$-bundle
$\tim\stackrel{\quo}{\longrightarrow} M$) if $s$ is measurable and $\quo\circ s={\rm id}_M$. 
Given a random variable $X:\Omega\to M$,  we call a map $\tX:\Omega\to \tim$ a {\em measurable lift} of $X$ 
if $\tX=s\circ X$ for some measurable section $s$.  
\end{defn}

Observe that a measurable lift $\tX$ of a measurable map $X:\Omega\to M$ is automatically 
measurable (see the discussion near the end of Section \ref{sect:intro}), hence is an $\tim$-valued random variable for which $\quo\circ \tX=X.$

The (left) principal $G$-bundle $\tim\stackrel{\quo}{\longrightarrow} M$ need not 
be trivial, hence may have no {\em continuous} sections  (continuous
maps $s:M\to\tim$ such that $\quo\circ s={\rm id}_M$).  
However, there always exist {\em measurable} sections:

\begin{lemma} 
\label{lem:meas_triv}
There exists a 
measurable section $s:M\to \tim$.  In fact, for any $p\in M$, there exists a measurable section that is smooth
on the complement of the cut-locus of $p$.
\end{lemma}

\pf Fix any $p\in M$ and any $\tilde{p}\in \quo^{-1}(p)$.  
As mentioned in Section \ref{sec:bary-prebary}, (i) the cut-locus $\CL:=\CL(p)$
is a closed set of volume zero,  (ii)  the set $\Dtan(p)\subset T_pM$ is star-shaped
with respect to $0_p$ (so its diffeomorphic image $\D:=\D(p)=M\backslash\CL$ is contractible),
and (iii) $M=\D\,\disjoint\, \CL$. Since $\quo:\tim\to M$ is a submersion, Lemma \ref{submersions_and_meas_zero} 
implies that the set $\tilde{\CL}:=\quo^{-1}(\CL)\subset\tim$ has volume zero.
Since $\D$ is contractible, the restricted bundle $\tim|_{\D}$ is trivial, so
there exists a continuous map $s_1:\D\to \tim$ such that $\quo\circ s_1={\rm id}_\D$.
Since $\quo:\tim\to M$ is a local diffeomorphism, the continuity of $s_1$ implies that
$s_1$ is smooth.  For each $q\in \CL$, arbitrarily select 
 a point $\tilde{q}\in \quo^{-1}(q)$, thereby constructing a map $s_2:\CL\to \tim$ satisfying 
$\quo\circ s_2={\rm id}_\CL$.  Then the map $s:M\to \tim$ defined by setting $s|_{\D}=s_1$
and $s|_{\CL}=s_2$ satisfies $\quo\circ s={\rm id}_M$, is continuous on $\D$, and is easily seen to be measurable.
\qedns

Lemma \ref{lem:meas_triv} can be generalized to any principal bundle over $M$; thus any such bundle is
{\em measurably trivial}. For the principal $G$-bundle 
$\tim\stackrel{\quo}{\longrightarrow} M$, this is essentially the content of Corollary
\ref{cor:meas_triv}(a) below.

\begin{cor} \label{cor:meas_triv} (a) Let $p\in M$. There exists an equivariant family 
$\{s_h\}_{h\in G}$ of measurable sections of $\tim\stackrel{\quo}{\longrightarrow} M,$ each of which is smooth on the 
open,  generic set $\D(p)$. 
(b) Let $X$ be an $M$-valued a random variable.  Then an equivariant family $\{\tX^{(h)}\}_{h\in G}$ of measurable 
lifts of $X$ exists.  If 
$X$ is absolutely continuous, then 
an equivariant family 
of {\em absolutely continuous} measurable lifts exists.

\end{cor}

\pf  (a)  Let $s:M\to \tim$ be a measurable section that is continuous on $\D(p)$.  
For each $h\in G$ let $s_h=h\dotprod s$. Then $\{s_h\}_{h\in G}$ is an equivariant family of measurable sections of $\tim\stackrel{\quo}{\longrightarrow} M,$ each of which is continuous on $\D(p)$.  

(b) Select any $p\in M$ and let $\{s_h\}_{h\in G}$ be as above.  
For each $h\in G$ let $\tX^{(h)}=s_h\circ X.$ Then 
$\{\tX^{(h)}\}_{h\in G}$ is an equivariant family of measurable lifts of $X$. 

Assume now that $X$ is absolutely continuous, and let $\mu_X$ be
the probability measure on $M$ induced by $X$. Let $E\subset \tim$
be a set of volume zero, and let $h\in G$.  Then $(\tX^{(h)})^{-1}(E)=X^{-1}(s_h^{-1}(E))
=X^{-1}(\quo(E)),$  so $\mu_{\tX^{(h)}}(E)=\mu_X(\quo(E)).$ Since the map $\quo$ is a local diffeomorphism, $\quo(E)$ has volume zero in $M$, implying that $\mu_X(\quo(E))=0$ (since $\mu_X$ is absolutely continuous).  Hence $\mu_{\tX^{(h)}}$ is absolutely continuous.
\qedns

\begin{remark}
\label{rem:assocrv}\rm
Suppose that $X$ is an absolutely continuous $M$-valued random variable
whose induced measure $\mu_X$ is supported in an open ball $B=B_r(p)$, where $r<\rinj(M,g)$.  
Let $\{\tilde{B}_h\}_{h\in G}$ be the even covering of $B$
given by Corollary \ref{evenlycov}, relabeled as in Remark \ref{relabel}.
Let $\{s_h\}_{h\in G}$ be an equivariant family of measurable sections of 
$\tim\stackrel{\quo}{\longrightarrow} M$ for which $s_h|_{\D(p)}$ is smooth (for each $h\in G$).  Note that $B\subset \D(p)$ since $r<\rinj(M,g)$.  Then by equivariance, 
each local section $s_h$ must be precisely the map $(\quo|_{\tilde{B}_h})^{-1}$ (viewed as
a map from $B$ to $\tim$).  Furthermore, the measure
$\mu_{\tX^{(h)}}$ is supported in $\tilde{B}_h$ and is canonically determined by $X$, in the sense that
$\mu_{\tX^{(h)}}$ depends only on the canonically defined local section $s_h|_B$ (and on $X$).
In view of this observation, we make the following simplifying definition:
\end{remark}

\begin{defn}\label{def:assocrv}\rm  In the setting of Remark \ref{rem:assocrv}, we refer to any equivariant family $\{\tX^{(h)}\}_{h\in G}$ determined by an equivariant family of measurable sections that 
are continuous on $B$, as {\em ``the'' associated lifts of $X$} to $\tim$-valued random variables
(implicitly identifying any two random variables that coincide a.e.).
\end{defn}

\begin{cor}
\label{cor:measzero.5}
Let $\tilde{A}\subset\tim$ be a $G$-invariant set of volume zero in $\tim$, let $A=\quo(\tilde{A})$, and write $\tim_*=\tim\setminus\tilde{A}$ and $M_*=\quo(\tim_*)=M\setminus A.$
Let $X, X_1, X_2,\dots, X_N$ be i.i.d. $M$-valued random variables
with induced probability measure $\mu_X$
and for $1\leq i\leq N$ let $\{\tX_i^{(h)}\}_{h\in G}$ be an equivariant family of measurable lifts of $X$. (Note that for each $h$, the random variables $\tX_1^{(h)}, \dots, \tX_N^{(h)}$ are i.i.d.)

\ss (a) For all $h,h_1,h_2\in G$ we have
\bearray
h_1\dotprod \FM(\{\tX_1^{(h_2)},\dots,\tX_N^{(h_2)}\})
&=&\FM(\{\tX_1^{(h_1h_2)},\dots,\tX_N^{(h_1h_2)}\})
\eearray
and
\bearray\quo(\FM(\{\tX_1^{(h)},\dots,\tX_N^{(h)}\}))&=&
\FM(\{X_1,\dots,X_N\}).
\eearray
Furthermore, 
$\EFM(\{X_1,\dots,X_N\})=\quo^{-1}(\FM(\{X_1,\dots,X_N\}) ).$

\ss (b) Assume now that $X$ is absolutely continuous.  Then 

\begin{itemize}
 \item[(i)]
Almost surely, $\FM(\{X_1,\dots,X_N\})\subset M_*$ and, correspondingly,
\linebreak
$\EFM(\{X_1,\dots,X_N\})\subset\tim_*$\,.

\item[(ii)]
For every $h\in G$,
almost surely  $\FM(\{\tX_1^{(h)},\dots,\tX_N^{(h)}\})\subset\tim_*$.

\end{itemize}

\ss (c) 
Let $r\in (0,\rcx(M,g)]$, and assume that 
$X$ is absolutely continuous and supported in an open ball $B$ of radius less than $r$. Then
\begin{itemize}
\item[(i)]
Almost surely, the Fr\'{e}chet mean set of $\{X_1,\dots,X_N\}$ is a single point and lies in  
\linebreak
$M_*\cap B$.  Correspondingly, the equivariant Fr\'{e}chet 
means of $\{X_1,\dots,X_N\}$ are unique up to the action of $G$, and lie in 
$\tim_*\cap\quo^{-1}(B)$.

\item[(ii)] Let $\{\tilde{B}_h\}_{h\in G}$ be the even covering of $B$ given by Corollary \ref{evenlycov}, with the sets relabeled as in Remark \ref{relabel},
and assume that the equivariant family  $\{\tX_i^{(h)}\}_{h\in G}$ consists of the associated lifts of 
$X_i$\,, $1\leq i\leq N$.
For every $h\in G$,
almost surely the Fr\'{e}chet mean set of $\{\tX_1^{(h)},\dots,\tX_N^{(h)}\}$
is a single point and lies in $\tim_*\cap \tilde{B}_h$.

\end{itemize}
\end{cor}

\pf  (a) This follows from Proposition \ref{prop:uniqueuptoG}(a).  

(b) 
Since $\quo$ is a local diffeomorphism, the set $A$ has volume zero in $M$.  Statements (i) and (ii)
therefore follow from Corollary \ref{cor:measzero.3}(a).  

(c) This follows from part (b) and Proposition \ref{prop:equivar_means}.
\qedns

Theorem \ref{thm1.2},  which is simply Corollary \ref{cor:measzero.5} minus the statements that explicitly 
involve the lifts $\{X^{(h)}\}$, is now proven.  

\begin{remark}\label{cor:measzero.5-strat}\rm
If $(\tim,\Tilde{\mathcal{T}})$ is a weakly stratified manifold in sense of Remark \ref{weakstrat}, and
for any $h\in G$ and any stratum $\Tilde{T}\in \Tilde{\mathcal{T}}$,
the set $h\dotprod \Tilde{T}$ is also a stratum, then $M$ inherits a quotient
(weak) stratification  by declaring the strata of $M$ to be the images of the strata of $\tim$ under $\quo$.
In Corollary \ref{cor:measzero.5}  we can then take $\tilde{A}$ to be the
union of all strata of $\tim$ other than the top stratum.  Then Corollary \ref{cor:measzero.5}
applies with $\tim_*$ and $M_*$ replaced by by $\tim^{\rm top}$ and $M^{\rm top}$, respectively.
\end{remark}

\setcounter{equation}{0}

\section{An application to ``partial scaling-rotation means''}
\label{sec:PSRmeans}

In this section we apply results from earlier sections to  a problem 
concerning symmetric positive-definite (SPD) matrices  that was the initial motivation behind Theorem
\ref{thm1.2}: establishing a strong relation between  {\em scaling-rotation {\rm (SR)} means} of samples of SPD matrices,
and the more easily computed {\em partial scaling-rotation {\rm (PSR)} means} of such samples, notions introduced in %\citep{rjgs}.
\cite{rjgs}.  To describe these means and the application of Theorem \ref{thm1.2}, we
first review the ``scaling-rotation framework'' introduced in earlier work of the authors.

Let $\Sym^+(m)$ denote the space of (real) $m\times m$ SPD matrices.  In 
%\citep{Jung2015} 
\cite{Jung2015} the authors defined a 
 ``distance-function" $\dsr$ (not a true metric) on $\Sym^+(m),$ studied further in  
 %\citep{gjs2017,gjs2021}, 
 \cite{gjs2017,gjs2021}, that measures rotation of eigenvectors and scaling of eigenvalues. 
For these purposes, an element  $S\in \Sym^+(m)$ is represented as an eigendecomposition,  a pair $(U,D)\in \tim(m):=SO(m)\times \Diag^+(m)$ for which $\F(U,D):=UDU^T=S$. (Thus $\F^{-1}(S)$  is the set of eigendecompositions
of $S$.)

\begin{figure}[t]
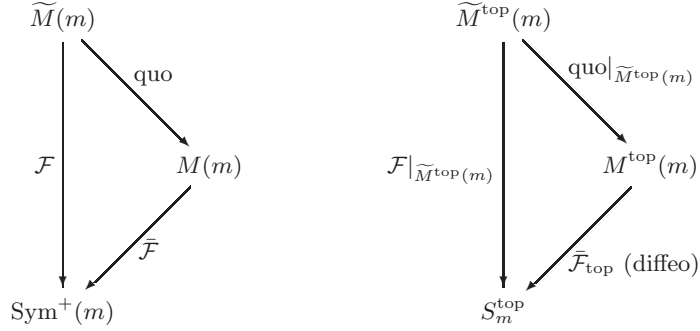

\begin{tabular}{cccc}
\mbox{\hspace{.4in}}
&
\mbox{
\begin{diagram}[nohug]
\tim(m) & &\\
&\rdTo>\quo&\\
\dTo^{\F}& &M(m) \\
&\ldTo>\bF&\\
\Sym^+(m)\\
\end{diagram}
}

&\mbox{\hspace{.5in}}&

\mbox{
\begin{diagram}[nohug]
\tim^{\rm top}(m) & &\\
&\rdTo >{\quo|_{\tim^{\rm top}(m)}}& \\
\dTo^{\F|_{\tim^{\rm top}(m)}}& &M^{\rm top}(m) \\
&\ldTo>{\bFtop\  (\mbox{diffeo})}&\\
\sptop
\end{diagram}
}
\end{tabular}

\caption{Commutative diagrams showing the relations among
various spaces and maps.}
\label{comdiag}
\end{figure}

The finite group $\G(m)$ of {\em even, signed permutations}, corresponding to 
orientation-preserving combinations of permutations and sign-changes of eigenvectors, 
acts freely on $\tim(m)$ (see \cite[Section 2.3]{gjs2017}), thereby determining a quotient manifold $M(m):=\tim(m)/\Sym^+(m)$ and a projection map $\quo:\tim(m)\to M(m)$.  
As depicted in the left-hand commutative diagram in Figure \ref{comdiag}, the manifold $M(m)$ sits ``between'' $\tim(m)$ and $\Sym^+(m)$: each $\G(m)$-orbit is contained in a fiber of the ``eigencomposition'' map $\F: \tim(m)\to \Sym^+(m), $ 
so $\F$ descends to a well-defined map $\bF:M(m)\to\Sym^+(m)$ satisfying
\be\label{Ffactors}
\F=\bF\circ\,\quo.
\ee
Note that if (and only if) $S\in \Sym^+(m)$ has distinct eigenvalues, then $\F^{-1}(S)$
is precisely an orbit of the $\G(m)$-action.  Thus the manifolds $\Sym^+(m)$ and $M(m)$ 
are not identical. However, they have full-volume-measure subsets that are ``essentially the
same''.  Specifically and more precisely, the subset $S_m^{\rm top}\subset\Sym^+(m)$ consisting of matrices with distinct eigenvalues, and the subset $\tim^{\rm top}(m):=\F^{-1}(\tim(m))\subset \tim(m),$ are open, generic
subsets of their ambient spaces. The set $\timtop(m)$ is $\G(m)$-invariant, and $S_m^{\rm top}$ 
can be identifed with the quotient of $\tim(m)^{\rm top}$ by the $\G(m)$-action. 
(We regard $\tim^{\rm top}(m), M^{\rm top}(m)$ as the ``top strata'' of their ambient spaces; see 
%\citep{gjs2017}.) 
\cite{gjs2017}.) 
Since $\tim^{\rm top}(m)$ is open in $\tim(m)$, its image $M^{\rm top}(m)$ under the covering-map $\quo$ is open in 
$M(m)$. The restriction of $\F$ to $\tim^{\rm top}(m)$ is a covering map $\tim^{\rm top}(m)\to \sptop,$ and 
the restriction of $\bF$ to $M^{\rm top}(m)$ is thus a diffeomorphism 
\ben
\bFtop:=\bF|_{M^{\rm top}(m)}:M^{\rm top}(m)\to \sptop\ ;
\een
see the right-hand commutative diagram in Figure \ref{comdiag}.  Thus $\bFtop$ naturally identifies
$M^{\rm top}(m)$ and $\sptop$ as manifolds, not just as point-sets.  It is in this sense that these subsets of 
the different spaces $M(m)$ and $\Sym^+(m)$ are ``essentially the same''.

The Lie groups $SO(m)$ and $\Diag^+(m)$ carry natural bi-invariant Riemannian metrics $g_{SO}$ and $g_{\D^+}$, respectively. (These are  complete and are canonical up to scale; for our choice of scale see \cite[Section 3.1]{gjs2017} .)
We define a product Riemannian metric $\tg$ on $\tim(m)$ by setting 
\lnbrk
$\tg=kg_{SO}\plus g_{\D^+,}$ where $k>0$ is an arbitrary relative-weight parameter.   
The action of $\G(m)$ on $\tim(m)$ is then isometric.
We are thus in the situation of Section \ref{sec:eqvt_means_general}: 
the manifold $M(m)$ then inherits a (unique) Riemannian metric $g$ for which the map 
$\quo:\tim(m)\to M(m)$ is a principal Riemannian covering map carrying $\tg$ to $g$.  Thus all the results of Section \ref{sec:eqvt_means_general} apply, with $\tim(m)$ as the covering space and $M(m)$  as the base space.
As in Section \ref{sec:eqvt_means_general},  we will denote the distance-functions (between points or subsets) arising from $\tg$ and $g$ as $\td$ and $d$, respectively.

Define functions $\dsr:\Sym^+(m)\times\Sym^+(m)\to\bfr$ (the {\em scaling-rotation distance})   and 
$\dpsr: \Sym^+(m) \times \tim(m)\to\bfr$ (the {\em partial scaling-rotation distance})   by
\bearray
\nonumber
\dsr(S_1,S_2) &=&\td(\F^{-1}(S_1),\F^{-1}(S_2)), \\
\dpsr(S, \tS_0) &=& \tilde{d}( \mathcal{F}^{-1}(S), \tS_0).
\label{defdpsr}
\eearray
Note that despite the formal similarity between equations \eqref{defdpsr} and \eqref{defdevt}, 
$\dpsr$ does not meet our definition of ``equivariant distance function''
 on its full domain, since $\tim(m)$ is a covering space only of $M(m)$, not of $\Sym^+(m)$.

\begin{remark}\label{rmk:isometry}
\rm  
For $S_1,S_2\in \sptop$,
\bearray\nonumber
\dsr(S_1,S_2)&=&\td(\F^{-1}(S_1),\F^{-1}(S_2))\\
\nonumber
&=& d(\quo(\F^{-1}(S_1)),\quo(\F^{-1}(S_2)))
\ \ \ \mbox{by Proposition \ref{quodist} }\\
&=& d(\bFtop^{-1}(S_1),\bFtop^{-1}(S_2)).
\label{dbareqdsr}
\eearray
This shows that $\dsrtop$, the restriction  of $\dsr$ to $\sptop\times\sptop$, is a metric
(precisely the result of 
\cite[Theorem 3.12] {Jung2015}, but proven by a somewhat different argument), 
and that, writing $\dtop$ for the restriction of $d$ to 
$M^{\rm top}(m)$,
\be\label{dbareqdsr.2}
\mbox{The diffeomorphism $\bFtop$
is an isometry
$(M^{\rm top}(m),\dtop)\to (\sptop, \dsrtop).$}
\ee

\end{remark}

For any configuration 
$\cals = (S_1,\ldots,S_N)\in (\Sym^+(m))^N$, we define functions \lnbrk 
$\fsrs:\sympp\to \bfr$ and $\fpsrs:\tim(m)\to\bfr$ by 
\bearray
\label{def_fsrs}
\fsrs(S)&=&\frac{1}{N}\sum_{i=1}^N  {\dsr}(S_i, S)^2, \\
\fpsrs(\tS)&=& \frac{1}{N}\sum_{i=1}^N  \dpsr( S_i, \tS)^2.
\label{def_fpsrs}
\eearray
\noi
As in \cite{rjgs}, 
we define the 
\emph{scaling-rotation} (SR) {\em mean set} and {\em partial scaling-rotation} (PSR) {\em  mean set}
of $\cals$, denoted $E_N^{\SR}(\cals)$ and $\tE_N^{\PSR}(\cals)$ respectively, by
\bestar\label{srmeanset}
E_N^{\SR}(\cals) &:=& 
\argmin(\fsrs)
\ \ \subset \Sym^+(m),\\
\label{psrmeanset}
\tE_N^{\PSR}(\cals) &:=& 
\argmin(\fpsrs)
\ \ \subset \tim(m).
\eestar
Corollary 4.10 of \cite{rjgs} establishes  that SR and PSR means always exist (i.e. that for any
$N$-point configuration $\cals$, the sets $E_N^{\SR}(\cals) $ and $\tE_N^{\PSR}(\cals)$ are nonempty). 

 Our main focus in this section will be on PSR means rather than SR means, for several reasons.
One reason is that the SR mean set of a general configuration $\cals$ can be very hard to compute, even when $\cals$ has a unique SR mean, since already the definition of $\dsr(\cdot,\cdot)$ involves a double-optimization.  As a practical alternative to using SR-means, 
%\citep{rjgs}  
\cite{rjgs}  
proposed using PSR  means to give more-easily computable approximations of SR means (upon 
projecting the PSR means 
to $\sympp$ via $\F$).   It is not clear, in general, how good this approximation is. However,
Theorem 3.5  of \cite{rjgs} established that in one important situation,
this approximation is perfect: 

\begin{thm}[{\cite[Theorem 3.5]{rjgs}}] \label{rjgsthm3.5}  Let $S_1, \dots S_N$ be elements of $\sptop$
and let $\cals=(S_1,\dots, S_N)$.  
If $E_N^{\SR}(\cals)\subset\sptop,$ then $\tE_N^{\PSR}(\cals) =\F^{-1}(E_N^{\SR}(\cals))$
and $\F(\tE_N^{\PSR}(\cals)) =E_N^{\SR}(\cals)$.
\end{thm} 

Since $S_m^{\rm top}$ 
 is an open, dense, full-measure subset of $\Sym( m),$ one might think that, 
  given i.i.d. absolutely continuous $\Sym^+(m)$-valued random variables $X_1,\dots, X_N$, 
  the genericity property ${\rm Pr}\big(E_N^\SR(X_1,\dots, X_N)\subset \sptop\big)=1$ would hold.
  This intuition is not correct,  essentially because the discontinuity of $\dsr$ leads to discontinuity  of the SR objective functions \eqref{def_fsrs}.  This is another reason we focus here on PSR means, which have nicer genericity properties.

Our final reason for focusing here on PSR means is that for a finite configuration that lies entirely in the top stratum 
$\sptop\subset \Sym^+(m)$,  the PSR mean set in $\tim(m)$ is exactly the equivariant Fr\'{e}chet mean set:

\begin{prop} \label{pmequivalence}
Let $S_1,\dots, S_n\in\sptop$, and for $1\leq i\leq N$, let $S_i^M=\bFtop^{-1}(S_i)$.  Let
$\cals=(S_1,\dots,S_N)$ and $\cals^M=(\cS_1,\dots,\cS_N)$.

Then 
\be\label{bfbsm}
\fpsrs =\tf_{\cals^M}
\ee
where the ``equivariant Fr\'{e}chet  function''
$\tf_{\cals^M}$ is defined as in \eqref{deftfq}. Hence $\tS$ is a PSR mean of $\cals$
if and only if $\tS$ is an equivariant Fr\'{e}chet mean of ${\cals}^M$, and
$\tE_N^\PSR(\cals)=\EFM(\cals^M)$.
\end{prop}

\pf
Using the definition of $\dpsr$, equation \eqref{Ffactors}, Proposition \ref{quodist}, the injectivity of $\bFtop$, and
equation \eqref{dpar=dbar}, we have
\bestar
\dpsr(S_i,\tS)=\td(\F^{-1}(S_i),\tS)
&=& \td(\quo^{-1}(\bF^{-1}(S_i)),\tS)\\
&=& d(\bF^{-1}(S_i),\quo(\tS))
\\
&=& d(\bFtop^{-1}(S_i),\quo(\tS))
\\
&=& d(\cS_i,\quo(\tS))
\\
&=& \devt(\cS_i ,\tS).
\eestar
Equation \eqref{bfbsm} now follows from the definition of $\tf_{\cals^M}$ (equation \eqref{deftfq}).
\qedns

Although the surjective maps $\bF:M(m)\to \sympp$ and $F:\tim(m)\to\sympp$ are not covering maps,
we can still define {\em measurable sections} for 
these maps, {\em measurable lifts} of $\Sym^+(m)$-valued random variables to $M(m)$-valued 
and $\tim(m)$-valued random variables, and 
($\G(m)$-){\em equivariant families of lifts}, to $\tim(m)$, of a $\sympp$-valued random variable:

\begin{enumerate}
\item Call a map $s:\sympp\to M(m)$ a {\em set-theoretic section} of $\bF$ if $\bF\circ s$ is the 
identity map of $\sympp$; equivalently, if $s(Y)\in \bF^{-1}(Y)$ for all $Y\in \sympp$. 
If, in addition, the map $s$ is 
measurable, call $s$ a 
{\em measurable section}.  
Similarly, call a map $\ts:\sympp\to\tim(m)$ a set-theoretic section or measurable
section of $\F$ if the appropriate condition(s) above hold with $(\bF,s)$ replaced by $(\F,\ts)$.

\ss
\item Call a map $\cX: \Omega\to M(m)$  (respectively, $\tX:\Omega\to\tim(m)$)
a {\em lift} of a map $X:\Omega\to\Sym^+(m)$ if $X=\bF\circ X^M$ (resp., $X=\F\circ\tX$); equivalently, 
if $X^M=s\circ X$ for some set-theoretic section $s$ of $\bF$
 (resp.,  if $\tX=\ts\circ X$ for some set-theoretic section $\ts$ of 
$\F$).  If $X$ is measurable, call $X^M$ (resp. $\tX^M$) a {\em measurable lift} of $X$ if 
$X^M=s\circ X$ (resp. $\tX^M=\ts\circ X$) for some measurable section $s$ of $\bF$ (resp. $\ts$ of $\F$). Note that measurable lifts are always measurable maps. 

\ss
\item Call an indexed collection of $\tim(m)$-valued random variables $\{\tX^{(h)}\}_{h\in \G(m)}$
an 
{\em equivariant family of lifts}  of a $\sympp$-valued random variable $X$ 
if (i) $\tX^{(h)}$ is a lift of $X$ for each $h\in \G(m)$, and (ii)
$\tX^{(h_1h_2)}=h_1\dotprod \tX^{(h_2)}$ for all $h_1,h_2\in \G(m)$.  

\end{enumerate}

Let 
$\splow=\Sym^+(m)\backslash \sptop$ and $\timlow(m)=\tim(m)\backslash \timtop(m)$.  
As detailed in \cite[Sections 2.4--2.7]{gjs2017}, the closed subsets 
$\splow$ and $\timlow(m)$  are finite unions of certain positive-codimension submanifolds (``lower strata'')
of $\Sym^+(m)$ and $\tim(m)$ respectively,  and hence have volume zero  in their ambient spaces.  Since $\quo$ 
is a local diffeomorphism, it follows that $M^{\rm low}(m):=\quo(\timlow(m))\subset M(m)$ is also a finite 
union of positive-codimension submanifolds and is a closed set.  Thus each of the sets $\splow, \timlow(m),$ 
and $M^{\rm low}(m)$ is a closed volume-zero subset of its ambient space. 

\begin{lemma}\label{liftlemma2} (a) Set-theoretic sections of $\bF$ exist, and all such sections are 
measurable. For any such section $s$, if $U\subset M(m)$ has volume zero, then so does $s^{-1}(U)$.

\ss 
(b)  Let $X$ be an absolutely continuous  $\Sym^+(m)$-valued
random variable $X$.  Then:

\begin{itemize}
\item[(i)]
Lifts of $X$ to $M(m)$ exist.  
All such lifts are absolutely continuous $M(m)$-valued random variables.  

\item[(ii)] An equivariant family $\{\tX^{(h)}\}_{h\in \G(m)}$ of measurable, absolutely continuous lifts of 
$X$ to $\tim(m)$ exists. 

\end{itemize}
\end{lemma}

For a proof of Lemma \ref{liftlemma2}, see Appendix \ref{sec:appendixA}.

\begin{prop} \label{prop:psr=efm}
Let
$X_1, \dots, X_N$ be i.i.d., absolutely continuous, $\Sym^+(m)$-valued
random variables, and let $\cX_1, \dots, \cX_N$ be lifts of $X_1, \dots X_N$ to $M(m).$
The following are true almost surely:

\begin{itemize}

\item[(a)]  All PSR means of $\{X_1,\dots,X_N\}$ lie in the top stratum:
\be\label{psrefm}
 \tE_N^\PSR(\{X_1,\dots, X_N\})=\EFM(\{X_1^M,\dots, X_N^M\}) \subset \timtop(m).
 \ee

\ss
\item[(b)] 
$\FM(\{\cX_1,\dots, \cX_N\})= 
\quo\big(\tE_N^{\PSR}(\{X_1,\dots, X_N\})\big)
 =\quo(\EFM(\{\cX_1,\dots, \cX_N\})$
\ \ and \newline
\ \ \ 
$\quo^{-1}\big(\FM(\{\cX_1,\dots, \cX_N\})\big)=
\tE_N^{\PSR}(\{X_1,\dots, X_N\}) =\EFM(\{\cX_1,\dots, \cX_N\})
$.

\end{itemize}
\end{prop}

\pf 
(a)  By Lemma \ref{liftlemma2}(b), the lifts $X_1^M, \dots, X_N^M$ are absolutely continuous $M(m)$-valued
random variables. For each Since $\timlow(m)=\quo^{-1}( M^{\rm low}(m))$ is 
 a $\G(m)$-invariant volume-zero subset of $\tim(m)$, 
 Corollary \ref{cor:measzero.5}(b) shows that the event (i) ``$\EFM(\{X_1^M,\dots, X_N^M\})\subset \timtop(m)$''
 occurs almost surely.  Since $\splow$ has volume zero and $X_1,\dots X_N$ are absolutely continuous,
the event (ii) ``$X_1,\dots, X_N$ all lie in $\sptop$'' occurs almost surely as well.
Hence the event ``(i) and (ii)'' also occurs almost surely.  But whenever (ii) occurs,  Proposition \ref{pmequivalence} 
implies that the equality in \eqref{psrefm} holds.   It follows that
 $\tE_N^\PSR(\{X_1,\dots, X_N\})\subset \timtop(m)$ almost surely.

\ss
(b) Proposition \ref{cor:measzero.5}(a) implies that we {\em always} (not just almost surely) have that
\lnbrk $\EFM(\{X_1^M,\dots, X_N^M\})=\quo^{-1}(\FM(X_1^M, \dots, X_N^M)),$ which, since $\quo$ is
surjective, 
implies that $\quo(\EFM(\{X_1^M,\dots, X_N^M\}))=\FM(X_1^M, \dots, X_N^M)$ as well.  Combining
this with the almost-sure statement \eqref{psrefm} yields the result.
\qedns

To apply results of Corollary \ref{cor:measzero.3} and Proposition \ref{prop:equivar_means} in the context of 
PSR means, we review some previously established geometric facts about the manifolds $(M(m),g)$ and $(\tim(m), \tg)$.  
As is well known, 
the manifold $(\Diag^+(m),g_{\D^+})$ has non-positive sectional curvature and infinite injectivity radius,
while $(SO(m),g_{SO})$ has non-negative sectional curvature, so
$\rinj(\tim(m),\tg)=\rinj(SO(m), kg_{SO})$ and $\Delta(\tim(m),\tg)=
\Delta(SO(m), kg_{SO})$. With the normalization chosen for $g_{SO}$ in \cite{rjgs}, 
$\rinj(SO(m), kg_{SO})=\sqrt{k}\pi $ and $\Delta(SO(m), kg_{SO})=k^{-1}/4.$
It is shown in \cite[Lemma 4.4]{rjgs} that $\inf\{\td(p,h\dotprod p) \ : \ 
h\in \G(m), p\in \tim(m)\}= \min\{d_{SO}(h,\id)\ : \  h\in \G(m)\setminus\{\id\}\}=:\bgp\leq  \frac{\pi}{2}$, so
\be\label{rcxmp}
\rcx(\tim(m),\tg)\ =\ \frac{\pi}{2} \sqrt{k}\ \geq \ \sqrt{k}\bgp.
\ee
Via Lemma \ref{lem:rinjbeta} and Corollary \ref{cor:rcxbeta},
the inequality in \eqref{rcxmp} then additionally yields
\be
\label{rinjmbar}
\rinj(M(m),g) 
\ =\  \sqrt{k}\bgp/2 \ \ \
\mbox{and}\ \ \ 
\rcx(M(m),g) 
\ =\  \sqrt{k}\bgp/4.
\ee

\begin{notation}\label{notn:9.4}\rm  (a) For $Y_0\in \sympp$ and $r>0$ let 
\bestar
\bar{\sf B}^{\SR}_r(Y_0) &=&\{Y\in\sympp : \dsr(Y,Y_0)\leq r\}, \\
{\sf B}^{\SR}_r(Y_0) &=&\{Y\in\sympp : \dsr(Y,Y_0)<r\}, \\
\mbox{and}\hspace{1in}
{\sf S}^{\SR}_r(Y_0) &=&\{Y\in\sympp : \dsr(Y,Y_0) = r\}.
\eestar
For $\tilde{Y}_0\in \tim(m)$ and $r>0$ let 
\ben
\tilde{\sf S}_r(\tilde{Y}_0; \td) = \{\tilde{Y}\in\tim(m) : \td(\tilde{Y},\tilde{Y}_0)= r\},
\een
the {\em sphere} with radius $r$ and center $\tilde{Y}_0$ in the 
Riemannian manifold $(\tim(m), \tg$). 

\ss (b) For $Y_0\in \sptop$ and $r>0$, we write $B_r(Y_0;\dsrtop)$ and $\bar{B}_r(Y_0;\dsrtop)$ for
 the open and closed balls, respectively, 
of radius $r$ and center $Y_0$ in the metric space $(\sptop,\dsrtop)$. 

\end{notation}

Informally, we 
call $\bar{{\sf B}}^{\SR}_r(Y_0),  {\sf B}^{\SR}_r(Y_0),$ and ${\sf S}^{\SR}_r(Y_0)$ the 
``closed $\dsr$-ball'',  ``open $\dsr$-ball'' and ``$\dsr$-sphere'' with center $Y_0$ and 
radius $r$, but we caution the reader not to let this terminology lead to implicit assumptions.  Since $\dsr$ is not a metric, 
these informally named balls and spheres are not balls and sphere in the metric-space sense, 
and cannot be expected to have all the properties of metric balls 
and spheres.  Indeed,  for $Y_0\in\sptop$, the function
$\dsr(\cdot, Y_0):\sympp\to\bfr$ is not even continuous 
(see \cite[Appendix A]{rjgs}).  It can be shown that 
$\bar{{\sf B}}^{\SR}_r(Y_0)$ is a closed {\em set}, at least, but 
${\sf B}^{\SR}_r(Y_0)$ is not always
an open set.  However, we will not need those 
facts to prove anything in the present paper.  Here, all that we need to know 
about ``$\dsr$-balls'' is in the following lemma.

\begin{lemma}\label{lem:balls}
Let $Y_0\in \sptop$ and let $r>0$. Then: 

\ss
(a) $\bar{B}_r(Y_0;\dsrtop)=\bar{\sf B}^{\SR}_r(Y_0)\cap \sptop$ and 
$B_r(Y_0;\dsrtop)={\sf B}^{\SR}_r(Y_0)\cap \sptop$.

\ss
(b) The set $\bar{{\sf B}}^{\SR}_r(Y_0)\minus
B_r(Y_0 ;\dsrtop)$ has volume zero in $\sympp$.

\ss
(c) Any absolutely continuous measure on $\Sym^+(m)$ supported in 
$\bar{{\sf B}}^{\SR}_r(Y_0)$ is also supported in 
$ B_r(Y_0;\dsrtop).$
\end{lemma}

\pf  (a) Immediate from definitions.

\ss 
(b) Let $\splow=\sympp\minus \sptop$. It is easily seen that 
\bestar
 \bar{{\sf B}}^{\SR}_r(Y_0)\minus B_r(Y_0,\dsrtop)
&\ \subset\ & \splow \cup {\sf S}_r^\SR(Y_0).
\eestar
Since $\splow\subset \sympp$ is a finite union of positive-codimension submanifolds, 
$\splow$ has volume zero in $\sympp$. Thus it suffices to show that 
${\sf S}_r^\SR(Y_0)$ also has volume zero.

\ss Let $\tilde{Y}_i, 1\leq i\leq |\G(m)|,$ be the elements of $\F^{-1}(Y_0)$. Then by
definition of $\dsr$, 
\be\label{union_spheres}
{\sf S}_r^\SR(Y_0) \subset \F\left( \bigcup_{i=1}^{|\G(m)|} 
S_r(\tilde{Y}_i ; \td)
\right).
\ee
 In any complete Riemannian manifold, the sphere with any given center and radius has volume zero in that manifold.
Hence each of the spheres $S_r(\tilde{Y}_i ; \td)$ has volume zero in $\tim(m)$,
and therefore so does their union.  Since $\dim(\tim(m))=\dim(\sympp)$ and $\F$ is smooth, 
$\F$ carries  volume-zero sets to volume-zero sets.  Hence ${\sf S}_r^\SR(Y_0)$
is contained in a volume-zero subset of $\sympp$, and therefore itself has volume zero in $\sympp$.

\ss (c) Immediate from part (a).
\qedns

\ss
In the next proposition, we consider  i.i.d., absolutely continuous, 
$\Sym^+(m)$-valued random variables $X_1, \dots, X_N$ whose induced measure $\mu_X$ is supported in 
$\bar{{\sf B}}_r^\SR(Y_0)$ for some given $Y_0\in\sptop$ and radius $r$ less than $\rcx(M(m),g)=\frac{\bgp\sqrt{k}}{4}$.
Lemma \ref{lem:balls}(b) assures us that $\mu_X$ is also supported in 
the open  metric ball $B_r(Y_0;\dsrtop).$ By Lemma \ref{liftlemma2}, for each $i$,
lifts $X^M_i$ of $X_i$ to $M(m)$ exist and are absolutely continuous,  and since $\bFtop:
(M^{\rm top}(m),\dtop)\to(\sptop, \dsrtop)$ is an isometry, 
the common measure $\mu_{_{\cX}}$ is supported in 
$ B_{r}(\bFtop^{-1}(Y_0))\subset (M(m),d)$. Since 
$r<\rcx(M(m),g)<\rinj(M(m),g)$, and $\quo:\tim(m)\to M(m)$ is a 
Riemannian covering map, each random variable $X_i^M$ has  an equivariant family of associated lifts
to $\tim(m)$ (see Definition \ref{def:assocrv}).  We call such families $\{\tX_1^{(h)}\}_{h\in \G(m)}, \dots, 
\{\tX_N^{(h)}\}_{h\in \G(m)}$, {\em consistently indexed}
if, in the setting of Remark \ref{rem:assocrv} and Definition \ref{def:assocrv}, 
there is an equivariant family of measurable sections $\{s_h\}_{h\in G}$
that simultaneously determines each of the families $\{\tX_i^{(h)}\}_{h\in \G(m)}.$

\begin{prop}\label{main_SR-PSR_result, version2}
Let $Y_0\in\sptop$,  let $Y_0'=\bFtop^{-1}(Y_0),$ let $r\in (0,\frac{\bgp\sqrt{k}}{4})$, and let 
$X_1, \dots, X_N$ be i.i.d., absolutely continuous, $\Sym^+(m)$-valued random variables whose distribution 
$\mu_X$ is supported in $\bar{{\sf B}}_r^\SR(Y_0).$ Let $\cX_1, \dots, \cX_N$ be lifts of 
$X_1, \dots X_N$ to $M(m)$, and let $\{\tX_1^{(h)}\}_{h\in \G(m)}, \dots, 
\{\tX_N^{(h)}\}_{h\in \G(m)}$ be consistently indexed, equivariant families of associated lifts of 
$\cX_1, \dots, \cX_N$ to $\tim(m)$; note that, thanks to the consistent indexing, for each $h\in\G(m)$ the 
random variables $\tX_1^{(h)}, \dots, \tX_N^{(h)}$ are i.i.d. 

The following are true almost surely:
\begin{itemize}
\item[(a)] For each $h\in \G(m)$, the set $\FM(\{\tX_1^{(h)},\dots,\tX_N^{(h)}\})$
consists of a single point and lies in $\F^{-1}(B_r(Y_0;\dsrtop))=\quo^{-1}
\left( B_r(Y_0')\cap M^{\rm top}(m)\right)\subset \tim^{\rm top}(m).$

\ss
\item[(b)] $\FM(\{\cX_1,\dots, \cX_N\})$ consists of a single point lying
in $B_r(Y_0')\,\cap\, M^{\rm top}(m).$  In addition,
$\tE_N^{\PSR}(\{X_1,\dots, X_N\})$ is identical to $\quo^{-1}\big(\FM(\{\cX_1,\dots, \cX_N\}\big)$ and
consists of a single $\G(m)$-orbit lying in $\F^{-1}(B_r(Y_0;\dsrtop)) \subset \tim^{\rm top}(m)$.   
\end{itemize}

\end{prop}

\ss \pf 
Let us write $\tim=\tim(m),$ $M=M(m)$, $\timtop=\timtop(m),$ and $M^{\rm top}=M^{\rm top}(m).$
Letting $B^M:=B_r(Y_0')\subset (M,d)$ play the role of the ball $B$ in Corollary
\ref{evenlycov}, let $\{\tilde{B}_h\}_{h\in \G(m)}$ be as in Corollary
\ref{evenlycov},  labeled equivariantly (see Remark \ref{relabel}).   
For each $\w\in \Omega$ (recall from Section \ref{sect:intro} that $\Omega$ is
the domain of all random variables in this paper), let us write $Q(\w):=(X_1(\w), \dots, X_N(\w))$ and
$Q^M(\w):=(X_1^M(\w), \dots, X_N^M(\w))$.

Since $r<\rcx(M,g)$ and $\tim\setminus \tim^{\rm top}$ is a $\G(m)$-invariant subset of $\tim$ of volume zero,
Lemma \ref{lem:balls} and Corollary \ref{cor:measzero.5}(c) together imply that there is a full-measure subset 
$\Omega_1\subset \Omega$ such that for all  $\w\in \Omega_1$:

\begin{itemize}
\item[(i)] For each $i\in \{1,\dots, N\},$ $X_i(\w)\in B_r(Y_0;\dsrtop)$ (implying that $X_i^M(\w)=\bFtop^{-1}(X_i(\w)),$ 
since $\bFtop$ is one-to-one). Equivalently, $Q(\w)\in (B_r(Y_0;\dsrtop))^N$,
which implies that $Q^M(\w)\in (B_r(Y_0'))^N$, since $\bFtop: (M^{\rm top},d^{\rm top})\to (\sptop,\dsrtop)$
is an isometry.

\ss
\item[(ii)] $\FM(Q^M(\w))$ 
is a single point and lies in  
$M^{\rm top}\cap B^M$.

\ss
\item[(iii)] For each $h\in G$, the  set $\FM(\tQ^{(h)}(\w))$
is a single point and lies in $\tim^{\rm top}\cap \tilde{B}_h$.   

\end{itemize}
Conclusion (a) of the Proposition is immediate from (iii). 
Conclusion (b) then follows from Corollary \ref{cor:measzero.5}(a) and Proposition \ref{prop:psr=efm}.
\qedns

\begin{remark}\label{weaker}\rm By simply  foregoing the first conclusion 
of Proposition \ref{main_SR-PSR_result, version2} we obtain the following weaker, but
more easily stated, result:

\begin{quotation} \em \normalsize
Let $Y_0\in\sptop$, let $r\in (0,\frac{\bgp\sqrt{k}}{4})$,
and let $X_1, \dots, X_N$ be absolutely continuous, i.i.d., $\Sym^+(m)$-valued
random variables whose distribution $\mu_X$ is supported in $\bar{{\sf B}}_r^\SR(Y_0).$ Then, almost surely,  
the partial scaling-rotation means of $X_1,\dots, X_N$ are unique up to the action of $\G(m)$
and lie in the top stratum of $\tim(m)$.
\end{quotation}
\end{remark}

\begin{remark}\label{rem:balls}\rm It can be shown that $\dsr$ is continuous on $\sptop\times\sptop$.
Hence, given $S_0\in \sptop$, $r>0,$ and $S\in \Bar{B_r(S_0;\dsrtop)}$ (the closure of $B_r(S_0;\dsrtop)$
in $\sympp$), 
\ben
\dsr(S,S_0)\leq \sup
\left\{\dsr(S',S_0): S'\in B_r(S_0;\dsrtop)\right\}=r;
\een
i.e. $S\in \bar{{\sf B}}_r^\SR(Y_0)$.  Thus $\Bar{B_r(S_0;\dsrtop)}\subset \bar{{\sf B}}_r^\SR(Y_0).$
Hence the support-condition on  $\mu_X$ in  Proposition \ref{main_SR-PSR_result, version2} and Remark \ref{weaker}
is met if $\mu_X$ is supported in $\Bar{B_r(S_0;\dsrtop)}.$
\end{remark}

{
\setcounter{equation}{0}
\appendix 
\section{Proofs of some technical results}
\renewcommand{\theequation}{\thesection.\arabic{equation}}
\label{sec:appendixA}

\noi{\bf Proof of Lemma \ref{submersions_and_meas_zero}}.
Let $n>0$. Suppose $A'\subset \bfr^n$ has volume zero in $\bfr^n$.
By \cite[Lemma 6.2]{lee_mfds}, $\bfr\times A'$ has volume zero in $\bfr^{n+1}$.
By induction, $\bfr^k\times A'$ has volume zero in $\bfr^{n+k}$ for all $k>0$. Trivially
this is true for $k=0$ as well.

Now let $M$ and $N$ be manifolds of dimensions $k+n$ and $n$ respectively,
and let $H:M\to N$ be a submersion. Identify $\bfr^{k+n}$ with $\bfr^k\times \bfr^n$ the
usual way. Since $H$ is a submersion, $M$ and $N$ have atlases
${\cal A},{\cal B}$ with the following property: for every chart $(U,\varphi)\in {\cal A}$
there exists $(V,\psi)\in{\cal B}$ such that $U=H^{-1}(V)$ and such that
the map $\psi\circ H\circ\varphi^{-1}: \linebreak
(\varphi(U)\subset \bfr^k\times \bfr^n) \to \bfr^n$ is simply the natural projection $\pi:\bfr^k
\times \bfr^n\to \bfr^n$ (restricted to $\varphi(U)\subset \bfr^k\times \psi(V)\subset \bfr^k\times \bfr^n$).

With all data as in the previous paragraph, suppose $A\subset  N$ has volume zero in $N$.
Then $A':=\psi(A\intersect V)$ has volume zero in $\bfr^n$.  Observing that
$\pi\circ\varphi\,(H^{-1}(A)\intersect U)=\pi\circ\varphi\,(H^{-1}(A\intersect V))=
\psi\circ H\,(H^{-1}(A\intersect V))=A'$, we have $\varphi(H^{-1}(A)\intersect U)
\subset \bfr^k\times A'$.The first paragraph of this proof therefore shows that
$\varphi(H^{-1}(A)\intersect U)$ has volume zero in $\bfr^{k+n}$. Since ${\cal A}$
is an atlas of $M$, this proves that $H^{-1}(A)$ has volume zero in $M$.\qedns

\ssn {\bf Proof of Lemma \ref{lem:rinjbeta}}.
 Let $\trho=\rinj(\tim,\tg)$ and $\rho=\rinj(M,g)$. 
Since $\quo$ is a local isometry, it intertwines
the exponential maps of $\tim$ and $M$:
\be\label{intertwine}
\quo\circ \exp_\tp =\exp_{\quo(\tp)}\circ\, \quo_{*\tp}\ \ \ \ \mbox{for all $\tp\in \tim$}.
\ee
It follows that for any  $p\in M,\  \tp\in \quo^{-1}(p),$ and $r>0$ for which $\exp_p:  B_r(0)\subset T_pM\to M$ 
is a diffeomorphism onto  its image, so is $\exp_{\tp}:  B_r(0)\subset T_{\tp}\tim\to \tim.$ Hence $\trho\geq \rho.$

Next, suppose that $\rho< \min\{\trho, \frac{\b}{2}\}$.
For $v\in T\tim$, let $\tgam_\tv$ denote
the geodesic $t\mapsto \exp(t\tv)$; for $v\in TM$ we similarly define the
geodesic $\g_{v}$.  We write $L[\cdot]$ for the length of a curve in $\tim$ or $M$.
Since $\quo$ is a local isometry, curve-lifting from $M$ to $\tim$ preserves length
and carries geodesics to geodesics.

Let $\tp\in \tim$, $p=\quo(\tp)$, and 
let $\trho_\tp, \rho_p$ denote the local injectivity radii of $\tim,M$ at $\tp,p$ respectively.
Since $\rho<\trho$, we may select $\e>0$ and $p\in M$ such that $\rho_p<\rho+\e<\trho$.
Then $\rho_p$ is finite, so  there exists $q\in M$ such that $q\in \CL(p)$ and $d(p,q)=\rho_p$.
Let ${v}\in T_{p}M$ be such that $\g_{v}|_{[0,1]}$ is a minimal geodesic
from $p$ to $q$, and let $\d= \e/\|v\|$; thus $L[\g_{v}|_{[0,1+\d]}]=\rho_p+\e<\trho$.
Since $q\in \CL(p)$, the geodesic $\g_{v}$ does not minimize  beyond $t=1$.  Hence
there exists $w\in T_{p}M$ such that $\exp_{p}(w)=\g_{v}(1+\d)$ and $\|{w}\|<(1+\d)\|{v}\|.$

Let $\a:[0,2+\d]\to \tim$ be the loop defined by $\a(t)=\g_{v}(t)$ for $0\leq t\leq
1+\d$ and $\a(t)=\g_{w}(2+\d-t)$ for $1+\d\leq t\leq 2+\d$.
Let $\tp\in \quo^{-1}(p)$.  Then the lift $\tilde{\a}$ of $\a$ to $\tim$ starting at $\tp$
satisfies $\tilde{\a}(2+\d)=h\dotprod \tp$ for some (unique) $h\in G$.
Let $\tq=\tilde{\a}(1+\d)$. Since $L[\tilde{\a}|_{[0,1+\d]}]=L[{\a}|_{[0,1+\d]}]<\trho$,
the curve $\tilde{\a}|_{[0,1+\d]}$ is a minimal geodesic from $\tp$ to
$\tq$, while $\tilde{\a}|_{[1+\d,2+\d]}$ is a shorter-length path from $\tq$ to $h\dotprod \tp$.  It follows
that $h\neq e$, and that
\be\label{2beta}
\b\leq \td(\tp,h\dotprod \tp)\leq L[\tilde{\a}]=L[\a]
<2(1+\d)\|v\|<2\rho+4\e.
\ee
Since \eqref{2beta} holds for arbitrarily small $\e$, this contradicts the hypothesis that $\rho<\b/2$.
Hence $\rho\geq \min\{\trho, \frac{\b}{2}\}.$ 

Finally, suppose the strict equality $\rho> \min\{\trho, \frac{\b}{2}\}$ holds.  Since $\rho\leq\trho$, we must
have $\min\{\trho,\frac{\b}{2}\}=\frac{\b}{2}$, so $\trho>\frac{\b}{2}$ and $\rho>\frac{\b}{2}$.
For each $p\in M$, let $\b_p$ be the minimal distance between distinct elements of $\quo^{-1}(p)$, and let $\tp_1, \tp_2\in \quo^{-1}(p)$ be a pair of points achieving this distance $\b_p$. Let $\tgam:[0,1]\to\tim$ be a minimal geodesic from 
$\tp_1$ to $\tp_2$ and let $\g=\quo\circ\tgam$.  Then $\g$ is a closed geodesic loop starting and ending at $p$, and
$L[\g]=L[\tilde{\g}]=\b_p$.  Hence $\rho_p\leq \frac{1}{2}L[\g]=\frac{\b_p}{2}$.  Taking the infimum over all 
$p\in M$ yields
$\rho\leq \frac{\b}{2}$, a contradiction. Hence $\rho= \min\{\trho, \frac{\b}{2}\}.$  
\qedns

\ssn {\bf Proof of Lemma \ref{liftlemma2}}.
(a) Using the Axiom of Choice, for each $Y\in \Sym^+(m)$ select an element
$s(Y)\in \bF^{-1}(\{Y\})$, thereby defining a set-theoretic section $s:\Sym^+(m)\to M(m)$.
Note that since $\bFtop$ is a bijection, $s|_{\sptop}=\bFtop^{-1}$.

Let $\call_1, \call_2$ denote the sigma-algebras of Lebesgue-measurable
subsets of $\Sym^+(m)$, $M(m)$ respectively (see Section \ref{sect:intro}).

For any $U\subset M(m)$, trivially $U=\big(U\cap M^{\rm top}(m)\big) \bigcup U\cap M^{\rm low}(m)$, and 
\linebreak therefore $s^{-1}(U)=s^{-1}\big(U\cap M^{\rm top}(m)\big) \bigcup 
s^{-1}(U\cap M^{\rm low}(m))$.  
Since $s^{-1}(U\cap M^{\rm low}(m)) \subset S_m'$, which has volume zero, 
it follows that $s^{-1}(U\cap M^{\rm low}(m))$ has volume zero, and therefore lies in $\call_1$.  Hence,
given $U\in \call_2$, to show that $s^{-1}(U)\in \call_1$ it suffices to show
that $s^{-1}\big(U\cap M^{\rm top}(m)\big)\in \call_1$.

First suppose that $U\subset M(m)$ is a Borel set.  Since $M^{\rm top}(m)$
is open, $U\cap M^{\rm top}(m)$ is a Borel set, and therefore so is 
$s^{-1}(U\cap M^{\rm top}(m))=s|_{\sptop}^{-1}(U\cap M^{\rm top}(m))=
\linebreak
\bFtop(U\cap M^{\rm top}(m))$.  Thus $s^{-1}\big(U\cap M^{\rm top}(m)\big)\in \call_1$. 

Next, suppose that $U\subset M(m)$ has volume zero.  Then every subset of $U$ has volume zero and therefore lies
in $\call_2$. In particular this applies to $U\cap M^{\rm top}(m)$. Again we have 
$s^{-1}(U\cap M^{\rm top}(m))=\bFtop(U\cap M^{\rm top}(m))$, 
and since $\bFtop$ is a diffeomorphism, it follows that $\bFtop(U\cap M^{\rm top}(m))$
has volume zero, hence lies in $\call_1$.  Furthermore, $s^{-1}(U)$ is the union
of two volume-zero sets, hence has volume zero.

Thus, for every $U\subset M(m)$ that is either a Borel set or a volume-zero set, 
$s^{-1}(U)\in \call_1$.   Since every $U\in {\mathcal L}_2$ is the union of a Borel set and a volume-zero set,
it follows that $s^{-1}(U)\in \call_1$ whenever $U\in\call_2$.  Thus $s$ is measurable.

\ss (b) Trivially, a map $X^M:\Omega\to M(m)$ is a lift of $X$ if and only if $X^M=s\circ X$
for some set-theoretic (hence measurable, by part (a)) section $s$ of $\bF$. 
 Since such sections exist, so do lifts $X^M$ of $X$.

Let $X^M$ be a lift of $X$ to $M(m)$, and let $s$ be the (necessarily unique) set-theoretic section of $\bF$
for which $X^M=s\circ X$.  Since $s$ and $X$ are measurable, so is $X^M$.

\ss
\hspace{.2in}
(i)
Let $U\subset M(m)$ be a volume-zero subset.  Then $s^{-1}(U)$ has volume zero, by part (a).  Since
$X$ is absolutely continuous, it follows that 
\ben
0=\mu_X(s^{-1}(U))= {\rm Pr}\big(X^{-1}(s^{-1}(U))\big)= {\rm Pr}\big( (s\circ X)^{-1}(U)\big)
=\mu_{X^M}(U).
\een
Hence the $M(m)$-valued random variable $X^M$ is absolutely continuous. 

\ss
\hspace{.2in}
(ii)  
By Corollary \ref{cor:meas_triv}(b) and its proof, there exists an equivariant family 
$\{\tX^{(h)}\}_{h\in \G(m)}$  of absolutely continuous, measurable lifts of $X^M$ from $M(m)$ to $\tim(m)$. 
Hence the family $\{\tX^{(h)}\}_{h\in \G(m)}$ has all the desired properties. 
\qedns

\section{Parametric transversality and a proof for Remark \ref{rmk:pi_NM_2}}
\label{sec:appendixB}
Given manifolds $\N$ and $\M$, and a submanifold $A\subset \M,$ we say that a smooth map $\F:\N\to \M$ 
is {\em transverse to $A$} if, for every $p\in \F^{-1}(A)$,
\ben
\F_{*p}(T_p \N) + T_{\F(p)}A = T_{\F(p)}\M
\een
(i.e. at $\F(p)$, the tangent space of $\M$  is spanned by the tangent space of $A$ together with
the image of $T_pN$ under the derivative of $\F$). If this condition is met, then $\F^{-1}(A)$
is a submanifold of $\N$, of codimension equal to the codimension of $A$ in $\M$. Note that 
trivially, a {\em submersion} $\F:\N\to \M$ is transverse to every submanifold of $\M$.

\begin{thm}[$C^\infty$ Parametric Transversality Theorem
{\cite[Theorem 6.35]{lee_mfds}}] \label{PTT}

Let $\N,\M,S$ be smooth manifolds and let $A\subset \M$ be a smooth submanifold.  Suppose
that $\{\F_s:s\in S\}$ is a smooth family of smooth maps $\N\to \M$;  i.e. assume that the map
$\F:S\times \N\to \M$ defined by $\F(s,p)=\F_s(p)$ is smooth. If $\F$ is transverse
to $A$,  then $\{s\in S: \F_s \ \mbox{is transverse to $A$}\}$ is a generic subset of $S$;
i.e. {\em generically} the maps $\F_s$ are  transverse to $A$.
\end{thm}

We will use Theorem \ref{PTT} to prove the following:

\begin{prop} For each $N\geq 1$, the set 
$\pi_N(\M_2)$ is generic in $M^N$.
\end{prop}

\pf 
Since $\M_1$ is open in $M^N\times M$, we can express $\M_1$ as a countable union of products
$U_i\times V_i$ ($i=1,2,\dots$), where $U_i$ and $V_i$ are open subsets of $M^N$ and $M$, respectively.
Since $F$ is a submersion, so is $F_i:=F|_{U_i\times V_i}$ for each $i$. Hence $F_i$ is
transverse to $Z$ for each $i$. For $Q\in U_i$, the evaluation-map $Y_{Q,i}: p\mapsto F_i(Q,p)$ is
simply $Y_Q|_{V_i}$. The proof of the Parametric Transversality Theorem 
as presented in \cite[Ch. 3, Theorem 2.7]{hirsch} shows that the set $A_i=\{Q\in U_i\ : \  Y_{Q,i} \ 
\mbox{is {\em not} transverse to}\ Z\}$ has volume zero in $U_i$, hence in $M^N$. Since there are only countably 
many $A_i$, their union has volume zero in $M^N$. But this union contains
$\{Q\in M^N\ : \  Y_Q \ \mbox{is not transverse to}\ Z\}$, and (using e.g. \eqref{YQ*p}) the condition 
``$Y_Q$ is transverse to $Z$'' is equivalent to ``$(\na Y_Q)|_p : T_pM \to T_p M$ is an isomorphism for 
all $p\in Z\intersect {\rm domain}(Y_Q)$''.
Hence the volume-zero set $\{Q\in M^N\ : \  Y_Q \ \mbox{is not transverse to}\ Z\}$ 
is exactly $M^N\setminus \pi_N(\M_2)$, and thus $\pi_N(\M_2)$ is generic in $M^N$.  
\qedns
}

\section*{Acknowledgements} The second author was supported by the National Research Foundation of Korea (NRF) 
grant funded by the Korea government (MSIT) (No. 2019R1A2C2002256).

\bibliographystyle{imsart-number}
%\bibliography{db}

%\section*{References}

\end{document}